%% file: The_Berger_conjecture_for_three_dimensional_manifolds.tex
\documentclass[a4paper,twoside]{article}

\usepackage{a4}
\usepackage[left=2.5cm,right=2.0cm,top=1.5cm,bottom=1cm,includeheadfoot]{geometry} 

\usepackage[ngerman, english]{babel}

\usepackage{lmodern}
\usepackage[T1]{fontenc}
\usepackage{textcomp}

\usepackage[utf8]{inputenc}

\usepackage{autobreak}
\usepackage{amsfonts}
\usepackage{amstext}
\usepackage{amssymb}
\usepackage{amsmath}
\usepackage{amsthm}
\usepackage{mathrsfs} 
\usepackage{wasysym}
\usepackage{mathtools}
\usepackage{appendix}

\usepackage{graphicx}
\usepackage{tikz}
\usepackage{tikz-cd}
\usepackage{tikz-3dplot}
\usetikzlibrary{matrix, calc, arrows,intersections,shapes.arrows}
\usepackage{pgfplots}
\pgfplotsset{
  compat=1.12,
  axis lines=middle
}

\usepackage{listings}
\usepackage{xcolor}

\usepackage{pdfpages}
\usepackage[pdfborder={0 0 0},colorlinks=false]{hyperref}
\usepackage[capitalize,nameinlink,noabbrev]{cleveref}

\usepackage[babel,german=quotes]{csquotes}
\usepackage[backend=biber, style=numeric]{biblatex}

\renewbibmacro*{in:}{}

\DeclareFieldFormat[article]{title}{\mkbibemph{#1}}
\DeclareFieldFormat[article]{journaltitle}{#1}
\DeclareFieldFormat[article]{pages}{#1}

\AtEveryBibitem{%
  \clearfield{year}%
  \clearfield{date}%
  \clearfield{isbn}%
 \clearfield{url}%
  \clearfield{urlyear}%
  \clearfield{doi}%
}

\renewbibmacro*{date}{}

\renewbibmacro*{issue+date}{}

\DeclareFieldFormat[article]{volume}{\mkbibbold{#1}}
\DeclareFieldFormat[article]{number}{\mkbibbold{#1}}

\input{theorems}

\usepackage{enumitem}

\newcommand{\comment}[1]{}
\input{commands}

\begin{document}

\input{abstract}

\tableofcontents

\input{introduction}


\input{topquotient}


\input{constructionhomologyclass}


\input{maintheoremproof}











\newpage

\printbibliography


\end{document}

%% file: theorems.tex
\makeatletter
\newtheoremstyle{theorem}
  {}
  {}
  {\itshape}
  {}
  {\bfseries}
  {.}
  { }
  {\thm@notefont{\bfseries}\thmhead@plain{#1}{#2}{#3}}
  
\newtheoremstyle{definition}
  {}
  {}
  {}
  {}
  {\bfseries}
  {.}
  { }
  {\thm@notefont{\bfseries}\thmhead@plain{#1}{#2}{#3}}  
\makeatother

\theoremstyle{definition}
\newtheorem{definition}{Definition}[section]
\crefname{definition}{Definition}{Definitions}

\crefname{assumption}{Assumption}{Assumptions}

\newtheorem{remark}[definition]{Remark}

\theoremstyle{theorem}

\crefname{theorem}{Theorem}{Theorems}
\newtheorem{theoremx}{Theorem}
 
\newtheorem{prop}[definition]{Proposition}
\crefname{prop}{Proposition}{Propositions}
\newtheorem{lemma}[definition]{Lemma}
\crefname{lemma}{Lemma}{Lemmas}
\newtheorem{corollary}[definition]{Corollary}
\crefname{corollary}{Corollary}{Corollaries}

%% file: commands.tex
\newcommand{\D}{\textup{d}}

\newcommand{\V}{\mathcal{V}}

\newcommand{\IZ}{\ensuremath{\mathbb{Z}}}
\newcommand{\IQ}{\ensuremath{\mathbb{Q}}}
\newcommand{\IR}{\ensuremath{\mathbb{R}}}
\newcommand{\IC}{\ensuremath{\mathbb{C}}}
\newcommand{\IS}{\ensuremath{\mathbb{S}}}

\DeclareMathOperator{\sgn}{sgn}
\DeclareMathOperator{\pr}{pr}
\DeclareMathOperator{\End}{End}

 {}  
\let \Im \undefined {} \DeclareMathOperator{\Im}{Im} 

\DeclareMathOperator{\ind}{ind} 			

\DeclareMathOperator{\stab}{Stab}

\DeclareMathOperator{\vspan}{span}

\let\originalleft\left
\let\originalright\right
\renewcommand{\left}{\mathopen{}\mathclose\bgroup\originalleft}
\renewcommand{\right}{\aftergroup\egroup\originalright}

%% file: abstract.tex
\title{The Berger conjecture for three dimensional manifolds}
\author{Jakob Dittmer}
\date{}

\maketitle

\begin{abstract}
A conjecture of Berger says that on a simply connected manifold all of whose geodesics are closed, all unit speed geodesics have the same shortest period. In this paper we will show that the $\IS^1$-action induced by the geodesic flow is free. We use symplectic geometry to investigate the topology of the quotient of the unit tangent bundle $T^1M/\IS^1$ by the $\IS^1$ action that comes from the geodesic flow. We show that this quotient admits the structure of a smooth symplectic manifold and that this symplectic manifold is symplectomorphic to $\IS^2\times \IS^2$ equipped with a split form. This will enable a computation involving the Euler class of the orbibundle $T^1M\rightarrow T^1M/\IS^1$ which will lead to a solution of Berger's conjecture for three-dimensional manifolds. 
\end{abstract}

%% file: introduction.tex
\section*{Introduction}

Given a Riemannian manifold $(M,\langle \cdot,\cdot\rangle)$, a geodesic $\gamma:I\rightarrow M$, defined on an interval $I\subset \IR$, is called a closed geodesic if it extends to a periodic map $\gamma:\IR\rightarrow M$. Riemannian manifolds all of whose geodesics are closed have long been objects of particular interest, especially since Zoll and Tannery discovered metrics on the sphere with this property that are different from the round metric. The most familiar examples are the compact symmetric spaces of rank $1$, the so-called CROSSes: spheres equipped with the round metric, real, complex, and quaternionic projective spaces equipped with their Fubini--Study metrics, and the Cayley plane.\\
A theorem of Wadsley states that if $(M,\langle \cdot,\cdot\rangle)$ is a Riemannian manifold all of whose geodesics are closed, then all geodesics have a common period. Consequently, the geodesic flow induces a smooth $\IS^1$-action on the unit tangent bundle $T^1M$. Berger's conjecture states that on a simply connected manifold with this property, all unit speed geodesics have the same minimal period. This is equivalent to the $\IS^1$-action induced by the geodesic flow on the unit tangent bundle being free. In this case, the volume of $M$ is an integral multiple of the volume of the round sphere whose radius is given by $2\pi$ divided by the common length of the geodesics. Gromoll and Grove proved the conjecture for two-dimensional manifolds in \cite{gromoll1981metrics}. In \cite{radeschi2017berger}, Radeschi and Wilking showed that it also holds for topological spheres of dimension at least $4$.\\
In dimension $3$, only partial results have been known so far. John Olsen proved that, for a simply connected three-dimensional manifold, the connected components of the isotropy set of the geodesic flow with fixed isotropy type are diffeomorphic to lens spaces whose fundamental groups have even order. In the earlier paper ``Isotropy sets of periodic geodesic flows on three-dimensional manifolds'', the author proved that, in this situation, the isotropy set has at least two connected components and, moreover, that all isotropy groups of the geodesic flow have odd order. The aim of the present paper is to prove Berger's conjecture also in dimension $3$. More precisely, we prove the following theorem.

\begin{theoremx}\label{bergconj}
Let $(M,\langle \cdot,\cdot\rangle)$ be a simply connected three dimensional Riemannian manifold all of whose geodesics are closed. Then all unit speed geodesics have the same shortest period.
\end{theoremx}

As already indicated in the abstract, our approach is based on a detailed study of the topology of the quotient of the unit tangent bundle by the $\IS^1$-action induced by the geodesic flow. In addition to Theorem A, we prove the following classification result.

\begin{theoremx}\label{quotclasstheorem}
Let $(M,\langle \cdot,\cdot\rangle)$ be a simply connected three dimensional manifold all of whose geodesics are closed. Then the quotient of the unit tangent bundle by the geodesic flow is homeomorphic to $\IS^2\times \IS^2$.
\end{theoremx}


Our strategy for proving Theorem A is as follows. Suppose that the isotropy set is nonempty, with connected components $C_1,\ldots,C_k$. The almost free $\IS^1$-action makes $T^1M/\IS^1$ a smooth orbifold with singular set $\bigcup_{i=1}^k C_i/\IS^1$. Each $C_i$ has isotropy group $\IZ/q_i\IZ$ with $q_i$ odd.\\

\noindent The symplectic structure on the cotangent bundle induces a symplectic orbifold form $\omega_{orb}$ that restricts to a symplectic form on each $C_i/\IS^1$. The de Rham theorem for orbifolds (see Satake) gives a class $[\omega_{orb}]$, and a theorem of Weinstein shows that $\operatorname{lcm}(q_1,\ldots,q_k)\cdot[\omega_{orb}]$ is integral. We compare two computations of the restriction of $[\omega_{orb}]$ to $C_i/\IS^1$.\\

\noindent On the one hand, $q_i[\omega_{orb}|_{C_i/\IS^1}]$ is the Euler class of the principal $\IS^1/(\IZ/q_i\IZ)$-bundle $C_i\rightarrow C_i/\IS^1$. Its evaluation on a suitably oriented fundamental class is the even integer $\#\pi_1(C_i)$.\\

\noindent On the other hand, a theorem of Weimin Chen on symplectic resolutions gives the quotient a smooth symplectic structure in which each $C_i/\IS^1$ is symplectic. The quotient is simply connected with second homology $\IZ^2$. We construct a class $c$ with $c^2\equiv 2\mod 4$ on which $[\omega_{orb}]$ vanishes, forcing the quotient to be homeomorphic to $\IS^2\times \IS^2$ and proving Theorem B.\\

\noindent We use Gromov's theory of pseudoholomorphic curves and the existence of at least two singular orbit types to construct a basis $E_1,E_2\in H_2(T^1M/\IS^1;\IZ)$ satisfying $E_1^2=E_2^2=0$ and $E_1\cdot E_2=1$, with $[C_i/\IS^1]=E_1+l_iE_2$ for $l_i\in\IZ$. For $k>2$, the numbers $l_i$ are all zero and this basis and the oddness of the $q_j$ show that $\operatorname{lcm}(q_1,\ldots,q_k)\cdot[q_i\omega_{orb}]$ has odd evaluation on $[C_i/\IS^1]$, contradicting the even evaluation above.\\

\noindent For $k=2$, a characteristic-class computation for the orbifold tangent bundle gives $\#\pi_1(C_1)+\#\pi_1(C_2)<4$, although both orders are positive and even. This completes the proof of Theorem A.\\


{\bf Acknowledgements}: This paper was to a large extent part of the author's doctoral thesis at the University of Münster. The author would like to express his deep gratitude to his advisor, Burkhard Wilking, for giving him the opportunity to work on this project, for his guidance and encouragement, and especially for suggesting how to attack the Berger conjecture. He would also like to warmly thank his mentor, Johannes Ebert, as well as Peter Albers and Michael Wiemeler, for many helpful conversations and valuable advice.

%% file: topquotient.tex
%


\section{Preliminaries}

\subsection{Geodesic flow}\label{sec::critmanhom}

Let $(M,\langle\cdot,\cdot\rangle)$ be a simply connected three-dimensional Riemannian manifold all of whose geodesics are closed. By Wadsley's theorem, there is a least common period $L>0$ of all unit-speed geodesics on $M$. Hence the geodesic flow induces an $\IS^1\cong \IR/L\IZ$-action on $T^1M$, and the Berger conjecture is equivalent to this action being free.\\
Since $M$ is three-dimensional, $M$ must be diffeomorphic to $\IS^3$, and $T^1M$ diffeomorphic to $\IS^3\times\IS^2$. By \cite{olsen2010three}, in this case the isotropy set of the geodesic flow is a finite disjoint union of finitely many connected compact three manifolds
\[
C_1\sqcup\cdots\sqcup C_k,
\]
Each $C_j$ is diffeomorphic to a lens space  has isotropy group $\IZ/q_j\IZ$, whose fundamental group has even order. Its quotient $C_j/\IS^1$ is diffeomorphic to $\IS^2$. Moreover by the author's prior paper the  isotropy set is not connected, if it is nonempty (See Theorem A in \cite{dittmer2026isotropysetsperiodicgeodesic}) so we may assume $k\geq2$ for the rest of the paper. Moreover, all $q_j$ are odd by Theorem B of the same paper. Since the action is almost free, the quotient $T^1M/\IS^1$ is a smooth orbifold.

\subsection{Free loop space}
The \emph{free loop space} $\Lambda M\coloneqq W^{1,2}(\IS^1,M)$ is a smooth Hilbert manifold. It carries the continuous circle action $([t],c)\mapsto(s\mapsto c(s+t))$. Evaluation at zero defines the fibration $p:\Lambda M\rightarrow M$, given by $p(c)=c(0)$. For further details, see \cite{klingenberg2012lectures}. Since $M\cong\IS^3$, we have
\begin{align}
H^*(\Lambda M;R)\cong H^*(\Lambda\IS^3;R)\cong H^*(\Omega\IS^3;R)\otimes H^*(\IS^3;R)
\end{align}
for any coefficient ring $R$. Here $\Omega\IS^3$ denotes the based loop space.

\subsection{Critical manifolds of the energy functional}
The energy functional is defined by
\[
E:\Lambda M\longrightarrow\IR,\qquad E(c)=\int_{\IS^1}|\dot c(t)|^2\,dt.
\]
It is smooth and $\IS^1$-invariant, and its critical points are precisely the closed geodesics \cite{klingenberg2012lectures}. It is also Morse--Bott, and we call the connected components of its critical set \emph{critical manifolds}. We denote the critical values by $0=e_0<e_1<\cdots$ and define
\[
K_i:=\operatorname{Crit}(E)\cap E^{-1}(e_i),
\qquad
\Lambda^i:=E^{-1}((-\infty,e_i+\epsilon_i]),
\]
where the numbers $\epsilon_i>0$ are chosen so that the intervals $(e_i-\epsilon_i,e_i+\epsilon_i)$ are pairwise disjoint. Each $K_i$ is connected by Corollary 3.3 of \cite{dittmer2026isotropysetsperiodicgeodesic}. The critical set $K_0$ consists of the constant loops and is therefore identified with $M$.\\

\noindent For a critical manifold $K$ of positive energy, consider the embedding
\[
\Phi:K\longrightarrow T^1M,\qquad
c\longmapsto\left(c(0),\frac{\dot c(0)}{\sqrt{E(c)}}\right).
\]
Its image is either all of $T^1M$ or one of the components $C_j$. We call $K$ \emph{regular} in the first case and \emph{singular} in the second. It is easy to see that a closed geodesic is regular precisely when its length is a multiple of $L$. If $K$ is singular, its lens-space structure gives
\begin{align}\label{eqrelativehomology}
H^*(K;\IZ)\cong
\begin{cases}
\IZ,& *=0,3,\\
\IZ/(\#\pi_1(K))\IZ,& *=2,\\
0,&\text{otherwise}
\end{cases}.
\end{align}
Here $\#\pi_1(K)$ is even. We write $c^q(t):=c(qt)$ for the $q$-fold iterate of $c$ and denote by $K^q$ the critical manifold consisting of the $q$-fold iterates of elements of $K$. Let $\widetilde K_j$ be the critical manifold of least energy whose image under $\Phi$ is $C_j$. It consists of prime geodesics, and its iterates give all critical manifolds mapping to $C_j$. The first regular critical manifold $K_{i_0}$ has energy $e_{i_0}=L^2$. The map $\Phi:K_{i_0}\rightarrow T^1M$ is an $\IS^1$-equivariant diffeomorphism.\\

\noindent For every critical manifold $B$ and every closed geodesic $c\in B$, the Hessian gives an orthogonal splitting
\begin{align}\label{eq:tangentspacedecomp}
T_c\Lambda M=T_cB\oplus T_c\Lambda M^-\oplus T_c\Lambda M^+.
\end{align}
The first summand is the kernel of the Hessian, while the second and third are its negative and positive spaces, respectively. The second space $T_c\Lambda M^+$ has always finite dimension which is called the \emph{index} of $c$. Since each critical manifold is connected, this index is constant along each critical manifold $B$ and the negative spaces form a finite-rank vector bundle $N(B)\to B$ with rank given by $\ind(B)$, called the \emph{negative bundle}. This bundle is in our case always orientable and the $\IS^1$-action preserves its orientation (See \cite{olsen2010three,radeschi2017berger} for more details about this). Morever, by Proposition 3.2 of \cite{dittmer2026isotropysetsperiodicgeodesic}, the index is strictly increasing with the energy, so geodesics with greater energy have also greater index.\\

\subsection{Symplectic geometry of the tangent bundle and the geodesic flow}
\label{sec::symplecticgeointro}

Choose an $\IS^1$-invariant metric $g$ on $T^1M$. For $p\in T^1M$, let $\V{\IS^1p}$ be the normal bundle of the $\IS^1$-orbit through $p$ with respect to the metric $g$. The stabilizer $\operatorname{Stab}(p)$ acts linearly and isometrically on the fiber $\mathcal{V}_p\IS^1p$. There exists $\varepsilon_p>0$ such that the exponential map
\begin{align*}
\exp:\V^{\varepsilon_p}\IS^1p\rightarrow T^1M
\end{align*}
restricts to an $\IS^1$-equivariant diffeomorphism onto its image. This induces a homeomorphism
\begin{align*}
\overline{\exp}:\V^{\varepsilon_p}_p\IS^1p/\operatorname{Stab}(p)\rightarrow \exp(\V^{\varepsilon_p}\IS^1p)/\IS^1.
\end{align*}
Thus $\overline{\exp}$ is an orbifold chart for the open neighborhood $\exp(\mathcal{V}^{\varepsilon_p}\IS^1p)/\IS^1$, and the collection of these maps forms an orbifold atlas of $T^1M/\IS^1$. Outside the set $\bigcup_{i=1}^k C_i/\IS^1$, where the $\IS^1$-action is free, these charts form a smooth manifold atlas of $T^1M/\IS^1\setminus \bigcup_{i=1}^k C_i/\IS^1$. 
\\

A symplectic orbifold form on $T^1M/\IS^1$ is induced by the canonical symplectic structure on the cotangent bundle. Let
\[
\pi_{M,*}:T^*M\longrightarrow M
\]
be the bundle projection. The Liouville form $\alpha\in\Omega^1(T^*M)$ is defined by
\[
\alpha_\eta(v)=\eta\bigl(D_\eta\pi_{M,*}(v)\bigr),
\qquad \eta\in T^*M,\; v\in T_\eta T^*M.
\]
In canonical coordinates $(x^1,\ldots,x^n,\xi^1,\ldots,\xi^n)$, one has
$\alpha=\sum_i\xi^i\,dx^i$; hence $-d\alpha$ is symplectic. The metric identifies $TM$ and $T^*M$ via
\[
\flat:TM\longrightarrow T^*M,\qquad
v\longmapsto \langle v,\cdot\rangle .
\]
Thus, for $\bar\alpha:=\flat^*\alpha$, the form $-d\bar\alpha$ is symplectic on $TM$. Let $Z$ be the geodesic vector field. Its restriction to $T^1M$ satisfies
$\bar\alpha(Z)=1$ and preserves $\bar\alpha$ (see \cite{besse2012manifolds}). Therefore, by Cartan's formula,
\[
0=\mathcal L_Z\bar\alpha
  =d\bigl(\bar\alpha(Z)\bigr)+\iota_Zd\bar\alpha
  =\iota_Zd\bar\alpha.
\]
Since $\mathcal L_Zd\bar\alpha=0$ as well, the form $d\bar\alpha|_{T^1M}$ is basic and descends to an orbifold $2$-form $\omega_{orb}$ on $T^1M/\IS^1$.\\

By Proposition 2.11 of \cite{besse2012manifolds}, for every $p\in T^1M$, the restriction of $d\bar\alpha_p$ to any complement of $\IR Z(p)$ in $T_pT^1M$ is nondegenerate. The transverse slices $\exp(\V_p\IS^1p)$ defining the orbifold charts provide such complements. Hence the descended form $\omega_{orb}$ is nondegenerate and therefore symplectic.

\section{Determing the homeomorphism type of the quotient \texorpdfstring{$T^1M/\IS^1$}{T1MS1}}
 
In this section, we determine the homeomorphism type of the quotient $T^1M/\IS^1$ under the effective $\IS^1$-action induced by the geodesic flow. We first show that, in addition to its natural structure as a symplectic orbifold, the quotient admits the structure of a symplectic manifold. We then determine its intersection form and use it to identify its homeomorphism type.

\subsection{The quotient \texorpdfstring{$T^1M/\IS^1$}{T1MS1} as a symplectic manifold}

\begin{prop}\label{propquotienthomeomorphismtyp}
The topological space $T^1M/\IS^1$ admits a smooth atlas $\mathcal{B}$ and a symplectic form $\omega\in \Omega^2(T^1M/\IS^1)$ such that
\begin{enumerate}[label=(\roman*)]
\item On $T^1M/\IS^1\setminus \bigcup_{i=1}^k C_i/\IS^1$, the form $\omega$ agrees with the orbifold $2$-form $\omega_{orb}$.
\item For each $i$, the form $\omega|_{C_i/\IS^1}$ agrees with $\omega_{orb}|_{C_i/\IS^1}$.
\end{enumerate}
In particular, each $C_i/\IS^1$ is a symplectic submanifold of $(T^1M/\IS^1,\mathcal{B})$.\\

\noindent Moreover, the quotient $T^1M/\IS^1$ is homeomorphic either to $\IS^2\times \IS^2$ or to $\IC P^2\#\overline{\IC P}^2$. If $T^1M/\IS^1$ admits a class $c\in H_2(T^1M/\IS^1;\IZ)$ with self-intersection $c^2\equiv 2\mod 4$, then $T^1M/\IS^1$ must be homeomorphic to $\IS^2\times \IS^2$.
\end{prop}

\begin{proof}
The quotient $T^1M/\IS^1$ is a four-dimensional symplectic orbifold whose singular set is the disjoint union of the symplectic two-spheres $C_1/\IS^1,\ldots,C_k/\IS^1$. By \cite{chen2018resolving}, such a symplectic orbifold admits a symplectic resolution with very good properties. More precisely, by Theorem 1 or Theorem 2 in \cite{chen2018resolving}, there is a smooth manifold atlas $\mathcal{B}$ on $T^1M/\IS^1$ and a symplectic form $\omega\in \Omega^2(T^1M/\IS^1)$ such that
\begin{enumerate}[label=(\roman*)]
\item On $T^1M/\IS^1\setminus \bigcup_{i=1}^k C_i/\IS^1$, the form $\omega$ agrees with the orbifold $2$-form $\omega_{orb}$,
\item For each $i$, the form $\omega|_{C_i/\IS^1}$ agrees with $\omega_{orb}|_{C_i/\IS^1}$.
\end{enumerate}
This proves the first part of the proposition. We now determine the homeomorphism type. Since $T^1M$ is diffeomorphic to $\IS^3\times \IS^2$, it is simply connected. Since the fibers of the quotient map are connected, \cite{calcut2012fundamental} implies that the quotient $T^1M/\IS^1$ is also simply connected. Moreover, because the $\IS^1$-action is almost free, the projection
\begin{align*}
Q:T^1M\times_{\IS^1} E\IS^1\rightarrow T^1M/\IS^1
\end{align*}
is a rational homotopy equivalence. It follows from Proposition 4.1 in "Geodesics flows on three dimensional manifolds with all geodesics closed" that $H^2(T^1M/\IS^1;\IQ)\cong \IQ^2$.
Since $(T^1M/\IS^1,\mathcal{B})$ is a compact, smooth, simply connected manifold, the second cohomology group $H^2(T^1M/\IS^1;\IZ)$ is a finitely generated free abelian group. Together with the preceding computation, this gives $$H^2(T^1M/\IS^1;\IZ)\cong \IZ^2.$$ In particular, the intersection form has rank $2$ and can only have signature $2$, $-2$, or $0$. Since $(T^1M/\IS^1,\mathcal{B})$ admits an almost-complex structure, a theorem of Wu (see 1.4.14 in \cite{gompf20234}) implies that the signature cannot be $2$ or $-2$. Hence the intersection form $Q_{T^1M/\IS^1}$ is indefinite. Freedman's classification of simply connected four-dimensional manifolds therefore implies that $T^1M/\IS^1$ is homeomorphic either to $\IS^2\times\IS^2$ or to $\IC P^2\#\overline{\IC P}^2$. In the latter case the intersection form is odd, and no class has square congruent to $2$ modulo $4$. Thus, if a class $c$ with $c^2\equiv 2\mod 4$ exists, the quotient must be homeomorphic to $\IS^2\times\IS^2$.
\end{proof}

%% file: constructionhomologyclass.tex
\subsection{Constructing a homology class with self-intersection \texorpdfstring{$\equiv 2 \mod 4$}{equiv 2 mod 4}}

We construct a class in $H_2(T^1M/\IS^1;\IZ)$ with self-intersection congruent to $2$ modulo $4$. For an open dense set of points $m\in M$, the tangent sphere $T^1_mM$ meets every component of the singular set transversally. For such $m$, we construct a smooth atlas $\hat{\mathcal{A}}$ on $T^1M/\IS^1$ for which the restricted orbit projection $\Pi_{\IS^1}:T^1_mM\rightarrow T^1M/\IS^1$ is an immersion and the antipodal map induces a smooth fixed-point-free involution. We then show that the normal Euler number satisfies $\chi(\V\Pi_{\IS^1})\equiv 2\mod 4$. Antipodal symmetry gives the same congruence for the self-intersection of $(\Pi_{\IS^1})_*([T^1_mM])$.

We begin with the transversality statement.

\begin{prop}\label{proptransgen} 
For generic $m\in M$ the unit sphere $T^1 _mM$ intersects each connected component $C\in \{C_1,...,C_k\}$ of the isotropy set transversally.
\end{prop}

\begin{proof}
Consider the canonical projection $\pi_M:T^1M\rightarrow M$, and fix $C\in\{C_1,\ldots,C_k\}$. By Proposition 3.6 of the authors prior paper "Isotropy sets of periodic geodesic flows on three-dimensional manifolds", the restriction $\pi_M|_C$ has non-zero mapping degree and is therefore surjective. By Sard's theorem, its regular values are dense in $M$; they also form an open set because $C$ is compact.\\

\noindent For any regular value $m$ and any $p\in C\cap\pi_M^{-1}(m)$, the differential $d(\pi_M|_C)_p$ is an isomorphism, since $C$ and $M$ are three-dimensional. Hence
\[
T_pC\cap T_p(T^1_mM)=T_pC\cap\ker d(\pi_M)_p=\{0\},
\]
so $T^1_mM$ intersects $C$ transversally. Taking the intersection of these open dense sets over the finitely many components $C_1,\ldots,C_k$ proves the claim.
\end{proof}

We first modify the orbifold atlas from Section 3.1.1 near the singular set to obtain a smooth manifold atlas. On the regular part of $T^1M/\IS^1$, we use exponential charts. Near each singular component $C$, we require the changes of orbifold charts to be complex linear in the second variable; see Section 2 of \cite{martin2021resolution} for the definition. For $m\in M$ satisfying the transversality condition, we also require the charts centered at $p\in T^1_mM\cap C$ to be Darboux charts for $\omega_{orb}$.


Fix a component $C$ of the singular set with isotropy group $\IZ/q\IZ$. Write
\begin{align*}
\V C =\{w\in TT^1M|_{C}:g(v,w)=0~\textup{for}~v\in T C,~ \pi_{TM}(v)=\pi_{TM}(w) \}
\end{align*} 
for the normal bundle of $C$ with respect to $g$, and
\begin{align*}
\V ({\IS^1p}\subset C)=\{ v\in TC:g(v,w)=0~\textup{for}~w\in T(\IS^1p),~ \pi_{TM}(v)=\pi_{TM}(w)  \}
\end{align*}
for the normal bundle of $\IS^1p$ in $C$. Since $C$ and $T^1M$ are orientable, so is $\V C$. Choose $\varepsilon>0$ so that the normal exponential map
\begin{align*}
\exp:\V ^{\varepsilon} C=\{v\in \V C:g(v,v)<\varepsilon ^2\}\rightarrow T^1M
\end{align*}
is a diffeomorphism onto its image and induces a homeomorphism
\begin{align*}
\bar{\exp}:\V ^{\varepsilon} C/\IS^1\rightarrow \exp(\V ^{\varepsilon}C)/\IS^1.
\end{align*} 
Similarly, choose $\varepsilon_p>0$ so that
\begin{align*}
\exp: \V^{\varepsilon_p}(\IS^1p\subset C)\coloneqq\{v\in \V{\IS^1p}:g(v,v)<\varepsilon_p ^2 \}\rightarrow C
\end{align*}
is a diffeomorphism onto its image. This induces a homeomorphism 
\begin{align*}
\bar{\exp}: \V^{\varepsilon_p} _{p}(\IS^1p\subset C)/(\IZ/q\IZ)=\V^{\varepsilon_p} (\IS^1p\subset C)/\IS^1 \rightarrow \exp(\V^{\varepsilon_p} (\IS^1p\subset C) )/\IS^1.
\end{align*}

For a smooth curve $\gamma:[0,1]\rightarrow T^1M$, let $\mathcal{P}_\gamma$ denote parallel translation with respect to $g$. As a component of the fixed-point set of an isometric action, $C$ is totally geodesic. 
Thus, for curves $\gamma$ contained in $C$,

\begin{align*}
\mathcal{P} _\gamma(T_{\gamma(0)}C)\subset T_{\gamma(1)}C
\end{align*}
and
\begin{align*}
\mathcal{P} _\gamma(\V_{\gamma(0)}C)\subset \V_{\gamma(1)}C.
\end{align*}
Consequently, the map
\begin{align*}
\exp(\V ^{\varepsilon_p} _p(\IS^1m\subset C))\times \V_pC &\rightarrow \V C|_{\exp(\V^{\varepsilon_p} _p(\IS^1p \subset C))},\\
(\exp(v'),v)&\mapsto \mathcal{P} _{t\mapsto \exp(tv')} (v) 
\end{align*}
is a smooth local trivialization of $\V C$.
The isotropy group $\IZ/q\IZ$ acts effectively by rotations in the second factor. The equivariant map
\begin{align}\label{eqorbifoldchartlift}
\phi_p:\V_p ^{\varepsilon_p}(\IS^1p\subset C) \times \V_p ^{\varepsilon} C&\rightarrow \exp\left(\V^{\varepsilon}C|_{\exp(\V_p ^{\delta_p}(\IS^1p\subset C))}\right)\nonumber\\
(v,w)&\mapsto \exp(\mathcal{P} _{t\mapsto \exp(tv)} (w))
\end{align}
induces a homeomorphism 
\begin{align}\label{eqorbifoldchart}
\tilde{\phi}_p:\V^{p} _p(\IS^1p\subset C)\times \V_p ^{\varepsilon} C/(\IZ/q\IZ)\rightarrow &\exp\left(\V^{\varepsilon}C|_{\exp(\V_p ^{\varepsilon_p}(\IS^1p\subset C))}\right)/(\IZ/q\IZ)\nonumber\\
&=\exp\left(\V^{\varepsilon}C|_{\exp(\V ^{\varepsilon_p}(\IS^1p\subset C))}\right)/\IS^1.
\end{align}
Thus $\tilde{\phi}_p$ is an orbifold chart of $T^1M/\IS^1$. To obtain manifold charts, use the map $z\mapsto z^q$, which induces a homeomorphism $\IC/(\IZ/q\IZ)\rightarrow\IC$ and restricts to a diffeomorphism
\begin{align*}
\IC\setminus \{0\}/(\IZ/q\IZ)\rightarrow \IC\setminus \{0\}.
\end{align*}
The maps
\begin{align}\label{eqmanifoldchart}
\bar{\phi}_p:\V^{\varepsilon_p} _p(\IS^1p\subset C)\times B_\varepsilon(0)&\rightarrow \exp\left(\V^{\varepsilon^{q}}|_{\exp(\V ^{\varepsilon_p}(\IS^1p\subset C))}C\right)/\IS^1,\nonumber\\
(v,w)&\mapsto \exp(v,(z\mapsto z^q) ^{-1}(w))
\end{align}
are therefore homeomorphisms onto open subsets of $T^1M/\IS^1$ and define manifold charts.

\begin{prop}\label{propfirstdesingquot}
Let $C$ be a component of the singular set. For each $p_0\in C$ there is $\delta_{p_0}>0$ such that if we consider $p,\bar{p}\in C$ with 
\begin{align}
\exp(\V^{\delta_p}(\IS^1p\subset C))\cap \exp(\V^{\delta_{\bar{p}}}(\IS^1\bar{p}\subset C))\neq \emptyset,
\end{align} 
the changes of the orbifold charts $\phi_p$ and $\phi_{\bar{p}}$ and the corresponding changes between the manifold charts $\bar{\phi}_p$ and $\bar{\phi}_{\bar{p}}$ are smooth. 
\end{prop}

\bigskip

\begin{proof}[Proof of Proposition \ref{propfirstdesingquot}]
Let $C\in \{C_1,\ldots,C_k\}$ be a component of the isotropy set of the geodesic flow. For sufficiently small $\delta_p>0$, it suffices to show that the changes between overlapping orbifold charts $\tilde{\phi}_p$ and $\tilde{\phi}_{\bar{p}}$ admit smooth lifts. Thus we seek open sets
\begin{align*}
V_{p,\bar{p}}&\subset \V_p ^{\delta_p}(\IS^1p\subset C) \times \V_p ^{\varepsilon} C,\\
V_{\bar{p},p}&\subset \V_{\bar{p}} ^{\delta_{\bar{p}}}(\IS^1\bar{p}\subset C) \times \V_{\bar{p}} ^{\varepsilon} C,
\end{align*}
and diffeomorphisms
\begin{align*}
\phi_{p,\bar{p}}:V_{p,\bar{p}}&\rightarrow V_{\bar{p},p},\\
\phi_{\bar{p},p}:V_{\bar{p},p}&\rightarrow V_{p,\bar{p}}
\end{align*}
such that the diagram
\begin{center}
\begin{tikzpicture}[commutative diagrams/every diagram]
\matrix[matrix of math nodes, name=m,column sep = 0.8cm, commutative diagrams/every cell] {
 V_{p,\bar{p}}& V_{\bar{p},p}\\
 V_{p,\bar{p}}/(\IZ/q\IZ)& V_{\bar{p},p}/(\IZ/q\IZ)\\
};
\path[commutative diagrams/.cd, every arrow, every label]
(m-1-1) edge node {$\phi_{p,\bar{p}} $} (m-1-2)
(m-2-1) edge node {$\tilde{\phi}_{\bar{p}}^{-1}\circ \tilde{\phi}_p$} (m-2-2)
(m-1-1) edge (m-2-1)
(m-1-2) edge  (m-2-2)
;
\end{tikzpicture}
\end{center}
commutes, and similarly for the reverse transition. We first construct smooth $\IS^1$-valued functions that move points between the transverse slices, and then compare the normal-bundle trivializations.\\

For each $p\in C$, choose $\delta_p>0$ such that
\begin{align*}
B_{4\delta_p}(p)\subset e^{2\pi i(-\frac{1}{2q},\frac{1}{2q})}\cdot \exp(\V_p ^{\varepsilon_p}(\IS^1p\subset C)).
\end{align*}
Then $B_{2\delta_p}(p')\subset B_{4\delta_p}(p)$ for every $p'\in B_{2\delta_p}(p)$. Suppose that $p,\bar{p}\in C$ satisfy
\begin{align}\label{chartintersection}
\exp(\V^{\delta_p}(\IS^1p\subset C))\cap \exp(\V^{\delta_{\bar{p}}}(\IS^1\bar{p}\subset C))\neq \emptyset.
\end{align}
There are $\IZ/q\IZ$-invariant open sets
\begin{align*}
W_{p,\bar{p}}&\subset \V_p ^{\delta_p}(\IS^1p\subset C),\\
W_{\bar{p},p}&\subset \V_{\bar{p}} ^{\delta_{\bar{p}}}(\IS^1\bar{p}\subset C)
\end{align*}
with
\begin{align*}
\IS^1\cdot\exp(W_{p,\bar{p}})
&=\IS^1\cdot\exp(W_{\bar{p},p})\\
&=\exp(\V^{\delta_p}(\IS^1p\subset C))\cap \exp(\V^{\delta_{\bar{p}}}(\IS^1\bar{p}\subset C)).
\end{align*}
The exponential maps identify both $W_{p,\bar{p}}/(\IZ/q\IZ)$ and $W_{\bar{p},p}/(\IZ/q\IZ)$ with the orbit space of this overlap.\\

We claim that there is a smooth function $g_{p,\bar{p}}:W_{p,\bar{p}}\rightarrow \IS^1$ such that
\begin{align*}
g_{p,\bar{p}}(x)\cdot\exp(x)\in \exp(W_{\bar{p},p})
\end{align*}
for every $x\in W_{p,\bar{p}}$. The corresponding change of base charts is then smooth and given by
\begin{align*}
\bar{\exp}^{-1}_{\bar{p}}\circ\bar{\exp}_p(x)
=\exp^{-1}\bigl(g_{p,\bar{p}}(x)\cdot\exp(x)\bigr).
\end{align*}
To construct $g_{p,\bar{p}}$ and $g_{\bar{p},p}$, choose $p'\in\exp(W_{p,\bar{p}})$ and assume, without loss of generality, that $\delta_{\bar{p}}\geq\delta_p$. There is $g\in\IS^1$ with $gp'\in\exp(W_{\bar{p},p})$. Hence $gp\in B_{\delta_{\bar{p}}}(gp')\subset B_{2\delta_{\bar{p}}}(\bar{p})$, and
\begin{align*}
gB_{2\delta_p}(p)=B_{2\delta_p}(gp)\subset B_{2\delta_{\bar{p}}}(gp)\subset B_{4\delta_{\bar{p}}}(\bar{p}).
\end{align*}
Since $g\exp(\V_p ^{\delta_p}(\IS^1p\subset C))\subset gB_{\delta_p}(p)=B_{\delta_p}(gp)$, we obtain
\begin{align*}
g\cdot\exp(\V_p ^{\delta_p}(\IS^1p\subset C))
&\subset B_{2\delta_{\bar{p}}}(gp)\subset B_{4\delta_{\bar{p}}}(\bar{p})\\
&\subset e^{2\pi i(-\frac{1}{2q},\frac{1}{2q})}\cdot\exp(\V_{\bar{p}} ^{\varepsilon_{\bar{p}}}(\IS^1\bar{p}\subset C)).
\end{align*}
The map
\begin{align*}
\left(-\frac{1}{2q},\frac{1}{2q}\right)\times\exp(\V_{\bar{p}} ^{\varepsilon_{\bar{p}}}(\IS^1\bar{p}\subset C))
&\rightarrow e^{2\pi i(-\frac{1}{2q},\frac{1}{2q})}\cdot\exp(\V_{\bar{p}} ^{\varepsilon_{\bar{p}}}(\IS^1\bar{p}\subset C)),\\
(s,v)&\mapsto e^{2\pi is}v
\end{align*}
is a diffeomorphism. Thus there is a smooth function
\(
h_{p,\bar{p}}:W_{p,\bar{p}}\rightarrow \left(-\frac{1}{2q},\frac{1}{2q}\right)
\)
with
\begin{align}\label{onesideinclusion}
ge^{2\pi i h_{p,\bar{p}}(x)}\cdot\exp(x)\in \exp(\V^{\varepsilon_{\bar{p}}}_{\bar{p}}(\IS^1\bar{p}\subset C)).
\end{align}
The required function is therefore
\begin{align*}
g_{p,\bar{p}}:W_{p,\bar{p}}\rightarrow\IS^1,\qquad x\mapsto e^{2\pi i h_{p,\bar{p}}(x)}\cdot g.
\end{align*}
The existence of the function $h_{p,\bar{p}}:W_{p,\bar{p}}\rightarrow \left(-\frac{1}{2q},\frac{1}{2q}\right)$ also implies with (\ref{onesideinclusion}) the inclusion
\begin{align*}
g^{-1}\cdot\exp(W_{\bar{p},p})\subset e^{2\pi i(-\frac{1}{2q},\frac{1}{2q})}\cdot\exp(\V_p ^{\varepsilon_p}(\IS^1p\subset C)).
\end{align*}
Applying the same argument to the diffeomorphism
\begin{align*}
\left(-\frac{1}{2q},\frac{1}{2q}\right)\times\exp(\V_p ^{\varepsilon_p}(\IS^1p\subset C))
&\rightarrow e^{2\pi i(-\frac{1}{2q},\frac{1}{2q})}\cdot\exp(\V_p ^{\varepsilon_p}(\IS^1p\subset C)),\\
(s,v)&\mapsto e^{2\pi is}v
\end{align*}
gives a smooth function
\begin{align*}
h_{\bar{p},p}:W_{\bar{p},p}\rightarrow\left(-\frac{1}{2q},\frac{1}{2q}\right)
\end{align*}
with
\begin{align*}
g^{-1}e^{2\pi i h_{\bar{p},p}(x)}\cdot\exp(x)\in \exp(\V^{\varepsilon_p}_p(\IS^1p\subset C)).
\end{align*}
This yields $g_{\bar{p},p}$ in the same way.

For the orbifold chart $\phi_p$ introduced above, the map
\begin{align*}
(x,v)\mapsto g_{p,\bar{p}}(x)\cdot\phi_p(x,v)
\end{align*}
is smooth and, in the normal-bundle trivializations, is a linear isometry in the second variable. Consequently,
\begin{align*}
\phi_{p,\bar{p}}:(x,v)\mapsto\phi_{\bar{p}}^{-1}\bigl(g_{p,\bar{p}}(x)\cdot\phi_p(x,v)\bigr)
\end{align*}
is smooth and a linear isometry in the second variable. Hence there is a smooth function
$$G_{p,\bar{p}}:W_{p,\bar{p}}\rightarrow\IS^1$$
such that
\begin{align}\label{orbifoldchangeofcharts}
\phi_{\bar{p}}^{-1}\bigl(g_{p,\bar{p}}(x)\cdot\phi_p(x,v)\bigr)
=\bigl(\bar{\exp}_{\bar{p}}^{-1}\circ\bar{\exp}_p(x),G_{p,\bar{p}}(x)\cdot v\bigr).
\end{align}
The reverse lift is obtained in the same way. This proves the first part.\\

For the manifold charts, formula (\ref{orbifoldchangeofcharts}) gives
\begin{align*}
\bar{\phi}_{\bar{p}}^{-1}\circ\bar{\phi}_p(x,v)
=\bigl(\bar{\exp}_{\bar{p}}^{-1}\circ\bar{\exp}_p(x),(G_{p,\bar{p}}(x))^q\cdot v\bigr),
\end{align*}
so these changes of charts are smooth as well.
\end{proof}

\begin{remark}\label{compatiblewithstandardorbstructure}
The charts $\phi_p$ are compatible with the orbifold structure on $T^1M/\IS^1$ given by the exponential maps. Suppose that $p,\bar{p}\in C$ and that the images of $\tilde{\phi}_p$ and the exponential chart $\bar{\exp}$ centered at $\bar{p}$ overlap. As in the proof of Proposition \ref{propfirstdesingquot}, we find, locally on this overlap, open subsets $W_{p,\bar{p}}$ and $W_{\bar{p},p}$ of the domains of $\phi_p$ and $\exp$, respectively, and smooth functions
$$g_{p,\bar{p}}:W_{p,\bar{p}}\rightarrow \IS^1,\qquad g_{\bar{p},p}:W_{\bar{p},p}\rightarrow \IS^1.$$
The maps
\begin{align*}
x&\mapsto \exp^{-1}(g_{p,\bar{p}}(x)\cdot\phi_p(x)),\\
x&\mapsto \phi_p^{-1}(g_{\bar{p},p}(x)\cdot\exp(x))
\end{align*}
are smooth local lifts of $\bar{\exp}^{-1}\circ\tilde{\phi}_p$ and $\tilde{\phi}_p^{-1}\circ\bar{\exp}$, respectively. Thus the two orbifold charts are compatible.
\end{remark}

\begin{corollary}\label{corsmoothstructure1}
The collection of the charts $\bar{\exp}_p$ for $p\in T^1M\setminus \bigcup_{i=1} ^k C_i$ and of the charts $\phi_p$ for $p\in \bigcup_{i=1} ^k C_i$ is a smooth orbifold atlas of the quotient $T^1M/\IS^1$ we will denote by $\mathcal{A}$. \\
Moreover, the collection of the charts $\bar{\exp}_p$ for $p\in T^1M\setminus \bigcup_{i=1} ^k C_i$ and of the charts $\bar{\phi}_p$ for $p\in \bigcup_{i=1} ^k C_i$ is a smooth manifold atlas of the quotient $T^1M/\IS^1$ we will denote by $\hat{\mathcal{A}}$.
\end{corollary}

\begin{proof}[Proof of Corollary \ref{corsmoothstructure1} ]
By Proposition \ref{propfirstdesingquot} the collection of the charts $\phi_p$ for $p\in C$ is a smooth orbifold atlas of $\exp(\V^{\varepsilon}C)/\IS^1$. 
On the smooth part of the orbifold $(\exp(\V^\varepsilon C)\setminus C)/\IS^1$ the maps $\phi_p$ are compatible with the smooth structure of $T^1M\setminus (\bigcup_{i=1} ^k C_i )/\IS^1$. So we obtain a smooth orbifold atlas on $T^1M/\IS^1$. \\

Moreover by Proposition \ref{propfirstdesingquot} the collection of the charts $\bar{\phi}_p$ for $p\in C$ is a smooth manifold atlas of $\exp(\V^{\varepsilon}C)/\IS^1$. Furthermore the restrictions
\begin{align*}
\bar{\phi}_p:\V^{\delta_p} _p(\IS^1p\subset C)\times (B_\varepsilon(0)\setminus \{0\})\rightarrow (\exp\left(\V^{\varepsilon^{q}}C|_{\exp(\V ^{\delta_p}(\IS^1p\subset C))}\right)\setminus C)/\IS^1
\end{align*}
are smooth with respect to the manifold structure $(T^1M\setminus \bigcup_{i=1} ^k C_i)/\IS^1$. So we get a smooth manifold atlas $\hat{\mathcal{A}}$ on the space $T^1M/\IS^1$.
\end{proof}

We will use the following identity for the antipodal map $A:T^1M\rightarrow T^1M$, $v\mapsto-v$.

\begin{prop}\label{proprelationinvolutionsymplecticform}
The antipodal map $A$ satisfies $A^{*}(d\bar{\alpha})=-d\bar{\alpha}$.
\end{prop}

\begin{proof}
It suffices to show that $A^{*}\bar{\alpha}=-\bar{\alpha}$.
Let $A:T^{*}M\rightarrow T^{*}M$, $p\mapsto-p$, also denote the
antipodal map on $T^{*}M$. By linearity of $\flat$, we have
$\flat\circ A=A\circ\flat$. Since $\bar{\alpha}=\flat^{*}\alpha$,
it remains to show that $A^{*}\alpha=-\alpha$. For $p\in T^{*}M$
and $v_p\in T_pT^{*}M$, the definition of $\alpha$ gives
\begin{align*}
(A^{*}\alpha)_p(v_p)
&=\alpha_{-p}(dA_p(v_p))\\
&=-p\bigl(d(\pi_{M,*}\circ A)_p(v_p)\bigr)\\
&=-p\bigl(d(\pi_{M,*})_p(v_p)\bigr)
=-\alpha_p(v_p),
\end{align*}
where we used $\pi_{M,*}\circ A=\pi_{M,*}$. This proves the claim.
\end{proof}

\begin{prop}\label{propspecialmetric}
Let $m\in M$, such that $T^1 _mM$ intersects each component of the singular set transversally. Then there is an $O(2)$-invariant metric $g$ on $T^1M$ with the following property:\\

For $m\in T^1M$ let $\V_p\IS^1p$ denote the normal bundle of the $\IS^1$-orbit of $p$ and for $s>0$ let $\V_p ^{s}\IS^1p$ denote the set of vectors in $\V\IS^1p$ with length less than $s$. 
Then for each $p\in T^1 _mM\cap \bigcup_{i=1} ^k C_i$ there is $\varepsilon_p>0$, such that 
\begin{align*}
\exp:\V_p ^{\varepsilon_p}\IS^1p\rightarrow T^1M
\end{align*}
is a diffeomorphism onto its image, such that $\exp_p ^{*}d\bar{\alpha}|_{V^{p,\varepsilon_p}}$ on the neighborhood $V^{p,\varepsilon_p}$ and $\exp_p ^{*}g$ are constant.
\end{prop}

\begin{proof}
Let $m\in M$ be as in the proposition, let $C$ be a component of the isotropy set with isotropy group $\IZ/q\IZ$, $q>1$, and let $p\in C\cap T^1_mM$. Choose an $\IS^1$-invariant auxiliary metric $h$ on $T^1M$. Then
\begin{align*}
T_pT^1M=\vspan\{Z(p)\}\oplus\V_p\IS^1p,
\qquad
\V_p\IS^1p=\V_p(\IS^1p\subset C)\oplus\V_pC.
\end{align*}
Identify each summand of $\V_p\IS^1p$ with $\IC$ by a linear isometry. Then the isotropy group $\IZ/q\IZ$ acts by rotations on the second factor. In a first step, we use the equivariant Darboux--Weinstein theorem to construct an equivariant Darboux chart near $p$. We will make $g$ Euclidean in these coordinates and then extend it to an $O(2)$-invariant metric on $T^1M$.\\

To obtain this chart, we first choose $h$ so that $d\bar{\alpha}_p$ has the standard form on $\V_p\IS^1p$. Start with any $O(2)$-invariant metric. Since $d\bar{\alpha}$ is nondegenerate on $\V_p\IS^1p$, choose a $d\bar{\alpha}$-compatible almost complex structure $J$ and replace $h$ on this space by the metric induced by $d\bar{\alpha}$ and $J$. This gives the auxiliary metric
\begin{align}\label{equauxmetric}
h|_{Z(p)} \oplus d\bar{\alpha}_p(\cdot,J\cdot ).
\end{align}
Then $J$ is a skew-symmetric isometry on $\V_p\IS^1p$. Unless otherwise stated, all metric-dependent terms below refer to the metric in (\ref{equauxmetric}).\\

The almost complex structure respects the decomposition
\begin{align*}
\V_p\IS^1p=\V_p(\IS^1p\subset C)\oplus\V_pC,
\qquad J_p=J^T_p\oplus J^\perp_p.
\end{align*}
Identifying these summands with $\IC$ by linear isometries, we may assume that
\begin{align*}
J_p(v,w)=(iv,iw).
\end{align*}
Thus, on $\V_p\IS^1p$, we have $d\bar{\alpha}_p(u,v)=h(J_pu,v)=-h(u,J_pv)$.
For sufficiently small $\varepsilon_p>0$, the map
\begin{align*}
\exp_p:B_{\varepsilon_p}(0)\times B_{\varepsilon_p}(0)&\rightarrow T^1M\\
(v,w)&\mapsto \exp_p((v,w))
\end{align*}
is a diffeomorphism onto its image $\exp_p(\V^{\varepsilon_p}_p\IS^1p)$.
The pullback $\exp_p^*d\bar{\alpha}$ is symplectic and satisfies
\begin{align}\label{equationpointdarboux}
(\exp_p^*d\bar{\alpha})_{(0,0)}((v_1,w_1),(v_2,w_2))
&=d\bar{\alpha}_p((v_1,w_1),(v_2,w_2))\nonumber\\
&=\langle iv_1,v_2\rangle_{\IR^2}+\langle iw_1,w_2\rangle_{\IR^2}.
\end{align}
The right-hand side defines a constant symplectic $2$-form invariant under rotations on the second factor. Since $(0,0)$ is fixed by this action, equation (\ref{equationpointdarboux}) allows us to apply the equivariant Darboux--Weinstein theorem (see Theorem 3.2 in \cite{DwivediHermanJeffreyvandenHurk2019}). We obtain a $\IZ/q\IZ$-equivariant diffeomorphism
\begin{align*}
\varphi_p:\V^{\varepsilon_p}_p(\IS^1p\subset C)\times\V^{\varepsilon_p}_pC
\rightarrow\V^{\varepsilon_p}_p(\IS^1p\subset C)\times\V^{\varepsilon_p}_pC
\end{align*}
with
\begin{align*}
(\varphi_p^*\exp_p^*d\bar{\alpha})_{(z_1,z_2)}((v_1,w_1),(v_2,w_2))
=\langle iv_1,v_2\rangle_{\IR^2}+\langle iw_1,w_2\rangle_{\IR^2}
\end{align*}
for $(z_1,z_2)\in\V^{\varepsilon_p}_p(\IS^1p\subset C)\times\V^{\varepsilon_p}_pC$ and $(v_1,w_1),(v_2,w_2)\in\V_p(\IS^1p\subset C)\times\V_pC$. In particular, $\varphi_p^{-1}\circ\exp_p^{-1}$ is a Darboux chart for $d\bar{\alpha}$ on $\exp_p(\V_p^{\varepsilon_p}\IS^1p)$.\\

To construct $g$, observe that the map
\begin{align*}
x_p:\left(-\frac{1}{2q},\frac{1}{2q}\right)\times\V^{\varepsilon_p}_p(\IS^1p\subset C)\times\V^{\varepsilon_p}_pC&\rightarrow T^1M\\
(t,v)&\mapsto e^{2\pi it}\exp_p\circ\varphi_p(v)
\end{align*}
is a diffeomorphism onto its image $U_p$. Choose $g$ on $U_p$ so that $x_p^*g$ is the Euclidean metric. For $z\in\IS^1$, the set $zU_p$ has coordinates $x_p^{-1}\circ z^{-1}$. The transverse coordinates are Darboux coordinates for $d\bar{\alpha}$. Define $g$ on $zU_p$ by ${z^{-1}}^*g$, so that it is Euclidean in these coordinates. Since $\varphi_p$ and $\exp_p$ are $\IZ/q\IZ$-equivariant, this is well defined and extends $g$ to an $\IS^1$-invariant metric on $\exp(\V^{\varepsilon_p}\IS^1p)$.

Similarly, define $g$ on $-U_p$ by $A^*g$, so that it is Euclidean in the coordinates $x_p^{-1}\circ A$. By Proposition \ref{proprelationinvolutionsymplecticform}, we have $A^*d\bar{\alpha}=-d\bar{\alpha}$. Hence $\varphi_p^{-1}\circ\exp_p^{-1}\circ A$ is a Darboux chart for $-d\bar{\alpha}|_{-\exp_p(\V^{\varepsilon_p}_p\IS^1p)}$, and the form is again constant in these coordinates. Repeating this construction for each pair $p,-p\in T^1_mM\cap\bigcup_{i=1}^kC_i$ gives an $O(2)$-invariant metric on an $O(2)$-invariant neighborhood of this set. Extending it to an $O(2)$-invariant metric on $T^1M$ yields the desired metric.
\end{proof}

\begin{remark}\label{remdarbouxchart}
Let $p\in T^1_mM\cap C$ be as in Proposition \ref{propspecialmetric}, and decompose
$$\V_p\IS^1p=\V_p(\IS^1p\subset C)\oplus\V_pC.$$
On a sufficiently small product neighborhood of the origin, $\exp_p^*g$ is constant, so constant vector fields are parallel. Hence
\begin{align*}
\exp_p((v,w))=\exp(\mathcal{P}_{t\mapsto \exp(tv)}(w))=\phi_p(v,w).
\end{align*}
Thus the maps $\phi_p$ in \eqref{eqorbifoldchartlift} are Darboux charts.
\end{remark}

\bigskip

We next point out that $A$ induces a smooth fixed-point-free involution on $T^1M/\IS^1$. For $v\in T^1M$ and $t\in[0,L]$, reversibility of the geodesic flow gives
\begin{align*}
\zeta^t (-v)=-\zeta^{L-t}(v).
\end{align*} 
Equivalently, the $\IS^1$-action satisfies 
\begin{align}\label{relationinvgeodesicflow}
z(-v)=-(z^{-1}v)
\end{align}
for $z\in\IS^1$ and $v\in T^1M$. Let 
\begin{align}\label{O2matrix}
B= \begin{pmatrix}
-1 &0\\
0& 1
\end{pmatrix}.
\end{align}
Every element of $O(2)$ is uniquely of the form $O'$ or $BO'$, where $O'\in SO(2)\cong\IS^1$. Since $BO'={O'}^{-1}B$, relation (\ref{relationinvgeodesicflow}) extends the geodesic-flow action to a smooth $O(2)$-action by
\begin{align}\label{definitingequationO2action}
(B\cdot O')\cdot v = A(O' \cdot v) 
\end{align}
for $O'\in SO(2)$. Relation (\ref{relationinvgeodesicflow}) also shows that $A$ descends to the quotient. If $A(v)=zv$, choose $w\in\IS^1$ with $w^2=z$. Then $A(wv)=w^{-1}A(v)=wv$, contradicting the fact that $A$ has no fixed points. Thus the induced involution $\bar{A}$ is fixed-point-free, and the diagram
\begin{equation}
\begin{tikzpicture}[commutative diagrams/every diagram]
\matrix[matrix of math nodes, name=m, commutative diagrams/every cell] {
T^1M& & T^1M  \\
 T^1M/\IS^1 & & T^1M/\IS^1 \\};
\path[commutative diagrams/.cd, every arrow, every label]
(m-1-1) edge node {$A$}(m-1-3)
edge node {$\Pi_{\IS^1}$}(m-2-1)
(m-1-3) edge node {$\Pi_{\IS^1}$} (m-2-3)
(m-2-1) edge node {$\bar{A} $ } (m-2-3);
\end{tikzpicture}
\end{equation}
commutes.

\begin{prop}\label{propsmoothnessinvolution}
Let $\bar{A}$ be the fix-point-free involution induced by the antipodal map $A$. Then $\bar{A}$ is smooth with respect to the smooth structure given by the atlas $\hat{\mathcal{A}}$.
\end{prop}

\begin{proof}
On the regular part $T^1M/\IS^1\setminus\bigcup_{i=1}^k C_i/\IS^1$, the smoothness of $\bar{A}$ is clear. It remains to prove smoothness near each $C/\IS^1$. Let $p\in C$ for some $C\in\{C_1,\ldots,C_k\}$, and let $\IZ/q\IZ$ be its isotropy group. We have the decomposition
\begin{align*}
T_pT^1M=\ker(d\bar{\alpha}_p)\oplus\V_p\IS^1p
=\vspan\{Z(p)\}\oplus\V_p\IS^1p.
\end{align*}
Since $g$ is $O(2)$-invariant, it is invariant under the stabilizer $\IZ/q\IZ$, and $A$ is an isometry. The map $A$ preserves $C$ and, by (\ref{relationinvgeodesicflow}), maps $\IS^1p$ to $\IS^1(-p)$. Hence $dA_p$ maps $\V_p(\IS^1p\subset C)$ to $\V_{-p}(\IS^1(-p)\subset C)$ and $\V_pC$ to $\V_{-p}C$. Thus $A$ acts on the normal bundle $\V C$ by $(x,w)\mapsto(-x,dA_x(w))$. The equivariant tubular neighborhood theorem identifies $\V^{\varepsilon}C$ with a neighborhood of $C$, and the exponential map is $A$-equivariant.

Choose the chart neighborhoods
\begin{align*}
V_p&=\exp\left(\V^\varepsilon C|_{\exp(\V_p^{\delta_p}(\IS^1p\subset C))}\right),\\
V_{-p}&=\exp\left(\V^\varepsilon C|_{\exp(\V_{-p}^{\delta_{-p}}(\IS^1(-p)\subset C))}\right)
\end{align*}
as in Proposition \ref{propfirstdesingquot}. After shrinking $\delta_p$ to $\delta'_p>0$, the map $A$ takes $V_p$ diffeomorphically onto an open subset of $V_{-p}$. The construction of the charts by parallel transport then gives
\begin{align*}
\phi_{-p}^{-1}\circ A\circ\phi_p:
\V_p^{\delta'_p}(\IS^1p\subset C)\times B_\varepsilon(0)
&\rightarrow\V_{-p}^{\delta_{-p}}(\IS^1(-p)\subset C)\times B_\varepsilon(0),\\
(x,v)&\mapsto(dA_px,\tilde{A}(x)(v))
\end{align*}
for a smooth map $\tilde{A}:\V_p^{\delta'_p}(\IS^1p\subset C)\rightarrow O(2)$. Since its domain is connected, $\det\tilde{A}$ is constant, and we distinguish two cases.

If $\tilde{A}$ takes values in $SO(2)$, the identification $SO(2)\cong\IS^1$ gives a smooth map $z:\V_p^{\delta'_p}(\IS^1p\subset C)\rightarrow\IS^1$ such that
\begin{align*}
\phi_{-p}^{-1}\circ A\circ\phi_p(x,v)
=(dA_px,z(x)\cdot v),
\end{align*}
where the dot denotes complex multiplication. Passing to the manifold charts by taking the $q$th power in the second coordinate gives
\begin{align*}
\bar{\phi}_{-p}^{-1}\circ\bar{A}\circ\bar{\phi}_p:
\V_p^{\delta'_p}(\IS^1p\subset C)\times B_{\varepsilon^q}(0)
&\rightarrow\V_{-p}^{\delta_{-p}}(\IS^1(-p)\subset C)\times B_{\varepsilon^q}(0),\\
(x,v)&\mapsto(dA_px,z(x)^q\cdot v),
\end{align*}
which is smooth.

If $\det\tilde{A}=-1$, then $B\tilde{A}$ takes values in $SO(2)$, where $B$ is the matrix in (\ref{O2matrix}). Under the usual identification of $\IR^2$ with $\IC$, this matrix acts by $v\mapsto-\bar{v}$. Absorbing the fixed minus sign into a smooth map $z:\V_p^{\delta'_p}(\IS^1p\subset C)\rightarrow\IS^1$, we can therefore write
\begin{align*}
\phi_{-p}^{-1}\circ A\circ\phi_p(x,v)
=(dA_px,\overline{z(x)\cdot v}).
\end{align*}
Since complex conjugation commutes with taking $q$th powers, the induced map in the manifold charts is
\begin{align*}
\bar{\phi}_{-p}^{-1}\circ\bar{A}\circ\bar{\phi}_p(x,v)
=(dA_px,\overline{z(x)}^{q}\cdot\bar{v}),
\end{align*}
which is also smooth. Thus $\bar{A}$ is smooth near every $C_i/\IS^1$, and the proof is complete.
\end{proof}

The restriction of the orbit projection to $T^1_mM$ has the following properties.

\begin{prop}\label{propsmoothimmersion}
Let $m$ be a point be such that $T^1 _{m}M$ intersects each component $C$ of the singular set transversally. The restriction of the orbit space projection $\Pi_{\IS^1}:T^1 _{m}M\rightarrow T^1M/\IS^1$ has the following properties:
\begin{enumerate}[label=\arabic*)]
\item It is a smooth immersion between $T^1 _{m}M$ equipped with its manifold structure and $T^1M/\IS^1$ equipped with its orbifold structure such that $\Pi_{\IS^1} ^{*}(\omega_{orb})=0$.
\item Let $C$ be a component of the singular set of $T^1M$ with isotropy $\IZ/q\IZ$. For $p\in  T^1 _mM\cap C$ there is a disc $D_p\subset \V_p ^{\varepsilon}C$, a chart $\theta_p: T^1 _mM\cap U_p\rightarrow D_p$ and a map $f_p:D_p\rightarrow \V ^{\delta_p}(\IS^1p\subset C)$ which is a diffeomorphism onto its image, such that 
\begin{align*}
\tilde{\phi}_p ^{-1}\circ \Pi_{\IS^1}\circ \theta^{-1} _p(x)=(f_p(x),x)
\end{align*}
and
\begin{align*}
\bar{\phi}_p ^{-1}\circ \Pi_{\IS^1}\circ \theta^{-1} _p(x)=(f_p(x),x^q).
\end{align*}
In particular $\Pi_{\IS^1}$ is a smooth immersion between $T^1 _{m}M$ equipped with its manifold structure and $T^1M/\IS^1$ equipped with the manifold structure given by the atlas $\hat{A}$.
\end{enumerate} 
\end{prop}

\begin{remark}
The map $x\mapsto(f_p(x),x)$ is a smooth lift of $\tilde{\phi}_p^{-1}\circ\Pi_{\IS^1}\circ\theta_p^{-1}$. Its image $\{(f_p(x),x):x\in D_p\}$ is Lagrangian for the standard symplectic form on $\IC^2$.
\end{remark}

\begin{proof}
Let $m\in M$ be such that $T^1 _mM$ intersects each component of the isotropy set transversally. We first show that
\begin{align*}
\Pi_{\IS^1}:T^1 _mM\rightarrow T^1M/\IS^1
\end{align*}
is a smooth Lagrangian immersion with respect to $\omega_{orb}$ on the orbifold $T^1M/\IS^1$. Choose $\varepsilon_p>0$ so small that
$$\exp:\V^{\varepsilon_p}\IS^1p\rightarrow T^1M$$
is a diffeomorphism onto its image $\exp(\V^{\delta_p}\IS^1p)$. Let $p\in T^1 _mM$, and choose $q$ smaller than each of the isotropy orders $q_1,...,q_k$. Then
\begin{align*}
\left(\frac{-1}{2q},\frac{1}{2q}\right)\times \V^{\varepsilon_p}_p\IS^1p&\rightarrow T^1M,~(x,v )\mapsto e^{2\pi i x}\cdot \exp(v)
\end{align*}
is a diffeomorphism onto its image. If $\stab_p$ is the isotropy group of $p$, then $(\exp,\V^{\delta_p} _p\IS^1p,\stab (p))$ is an orbifold chart of $\exp(\V^{\delta_p}\IS^1p)/\IS^1$, and
\begin{align*}
\left(\frac{-1}{2q},\frac{1}{2q}\right)\times \V^{\varepsilon_p}_p\IS^1p\rightarrow \V^{\varepsilon_p}_p\IS^1p,~(t,v)\mapsto v
\end{align*}
is a smooth local lift of $\bar{\exp} ^{-1} \circ \Pi_{\IS^1}\circ ((t,v)\mapsto e^{2\pi i t}\exp_p(v))$. At the origin, the differential of the chart $(t,v)\mapsto e^{2\pi i t}\exp_p(v))$ maps the tangent space of $(\frac{-1}{2q},\frac{1}{2q})$ to the span of $Z(p)$ and maps $\V _p\IS^1p$ to $\V_p \IS^1p$. 
The geodesic vector field is transverse to $T^1 _mM$: if $\gamma$ is a geodesic on $M$ starting at $m$ with initial vector $v\in T^1 M$, then $\zeta^t(v)=\dot{\gamma}(t)$ and $Z(v)=\frac{d}{dt}\zeta^t(v)|_{t=0}$. Since $d\pi_M(Z(v))=\dot{\gamma}(0)\neq 0$, $Z(v)$ does not lie in the kernel. Hence $d\pi_M(T_v T^1 _mM)\cap \vspan(Z(v))=\{0\}$, and therefore $\ker d(\Pi_{\IS^1})_{p}\cap T_pT^1 _mM =\{0\}$. Consequently, tangent vectors $v\in T_pT^1 _mM$ are determined by their images under
$$p_{\V_p\IS^1p}:T_pT^1 _mM\rightarrow \V_p \IS^1p$$
so this projection is injective. It follows that $\Pi_{\IS^1}:T^1 _{m}M\rightarrow T^1M/\IS^1$ has a smooth local lift with injective differential at $p$ and is a smooth orbifold immersion.\\
To show that $\Pi_{\IS^1}:T^1 _mM\rightarrow T^1M/\IS^1$ is Lagrangian with respect to $\omega_{orb}$, recall that $\omega_{orb}$ is induced by $-d\bar{\alpha}$, so $\Pi_{\IS^1} ^{*}(\omega_{orb})=-d\bar{\alpha}$. For $\Pi_{\IS^1}$ to be Lagrangian, it therefore suffices to show that
\begin{align}\label{equationT1Misvanished}
-d\bar{\alpha}|_{T^1 _mM}=0
\end{align}
holds. For this, it is enough to check that $\bar{\alpha}|_{T _mM}=0$. Since $\bar{\alpha}= \flat ^{*}\alpha$ and $\flat(T_mM)=T^{*} _mM$, this follows from the vanishing of $\alpha|_{T^{*} _mM}$, which holds by the definition of $\alpha$.\\ 

We now show that $\Pi_{\IS^1}:T^1 _mM\rightarrow T^1M$ is also a smooth immersion for the atlas $\hat{\mathcal{A}}$ constructed in Proposition \ref{propfirstdesingquot}. We prove that $\Pi_{\IS^1}$ is locally an embedding. Since $\hat{\mathcal{A}}$ agrees with the orbifold smooth structure outside the isotropy set, it suffices to consider $\Pi_{\IS^1}$ at points $\Pi^{-1} _{\IS^1}(\bigcup_{i=1} ^k C_i )$. Let $p\in C$ for some $C\in \{C_1 ,...,C_k  \}$. The map
\begin{align*}
\Theta_p:\left(-\frac{1}{2q},\frac{1}{2q}\right)\times \V_p ^{\delta_p}(\IS^1p\subset C)\times \V_p ^{\varepsilon}C&\rightarrow T^1M, \\
(t,v,w)&\mapsto e^{2\pi i t}\cdot \exp(\mathcal{P}_{s\mapsto \exp (sv)}(w))
\end{align*}
is a diffeomorphism onto its image $\mathcal{U}_p\subset T^1M$. 
Its restriction
\begin{align*}
\Theta_p:\{0\}\times \V_p ^{\delta_p}(\IS^1p\subset C)\times \V_p ^{\varepsilon}C\rightarrow \exp(\V^{\varepsilon}C|_{\exp(\V_p ^{\delta_p}(\IS^1p\subset C))})
\end{align*}
is the orbifold chart $\phi$. Decompose $T_pT^1M$ as
\begin{align*}
T_pT^1M=\vspan{Z(p)}\oplus \V_p {\IS^1p}=\vspan{Z(p)}\oplus \V_p {(\IS^1p\subset C)}\oplus \V_pC
\end{align*}
and identify $T_{(0,0,0)}((-\frac{1}{2q},\frac{1}{2q})\times \V_p ^{\delta_p}(\IS^1p\subset C)\times \V_p ^{\varepsilon}C)$ by
\begin{align*}
T_{(0,0,0)}(\left(-\frac{1}{2q},\frac{1}{2q}\right)&\times \V_p ^{\delta_p}(\IS^1p\subset C)\times \V_p ^{\varepsilon}C)\\
&=T_0\left(-\frac{1}{2q},\frac{1}{2q}\right)\oplus T_0\V_p ^{\delta_p}(\IS^1p\subset C)\oplus T_0\V_p ^{\varepsilon}C\\
&=\IR \oplus \V_p ^{}(\IS^1p\subset C)\oplus \V_pC ,
\end{align*}
using the standard tangent-space identifications for products and vector spaces. Under these identifications,
\begin{align*}
(D\Theta_p)_0(\IR\times \{0\}\times \{0\})&=\vspan (Z(p)),\\ (D\Theta_p)_{(0,0,0)}(\{0\}\times \V_p ^{}(\IS^1p\subset C)\times \{0\})&=\V_p ^{}(\IS^1p\subset C)
\end{align*}
and
\begin{align*}
(D\Theta_p)_{(0,0,0)}(\{0\}\times \{0\}\times \V_pC)=\V_pC.
\end{align*}
Now consider the orthogonal projections
\begin{align*}
P_{\V_p C}:T_pT^1M\rightarrow \V_pC,~~~~P_{\V_p (\IS^1p\subset C)}:T_pT^1M\rightarrow \V_p (\IS^1p\subset C)
\end{align*}
and
\begin{align*}
P_{\V_p ^{}\IS^1p}:T_pT^1M\rightarrow \V_p ^{}\IS^1p=\V_p ^{}(\IS^1p\subset C)\oplus \V_pC.
\end{align*}
Choose $J\in \End(TT^1M)$ such that
\begin{align*}
-d\bar{\alpha}(\cdot,\cdot)=g(\cdot,J \cdot).
\end{align*}
By our choice of $g$, the almost complex structure $J$ is multiplication by $i$.
Consider the subspace $L_p=P_{\V_p ^{}\IS^1p}(T_pT^1 _mM)$ of $T_pT^1M$. Since $T^1 _mM$ intersects $C$ transversally, the projection $P_{\V_pC}$ gives an isomorphism
\begin{align}\label{transversalintersectionlagrangian}
P_{\V_pC}(P_{\V_p {\IS^1p}}(T_pT^1 _mM))=P_{\V_pC}(T_pT^1 _mM)=\V_pC.
\end{align}
In particular, $P_{\V_p {\IS^1p}}(L_p)$ is two-dimensional. By (\ref{equationT1Misvanished}) and the fact that $Z(p)$ lies in the kernel of $-d\bar{\alpha}$, this subspace is Lagrangian with respect to $d\bar{\alpha}|_{\V_p {\IS^1p}}$. Together with (\ref{transversalintersectionlagrangian}), this gives
\begin{align*}
P_{\V_p ^{}(\IS^1p\subset C)}(P_{\V_p {\IS^1p}}(T_pT^1 _mM))=\V_p {(\IS^1p\subset C)},
\end{align*}
and hence the composition of the projection
\begin{align*}
P_{\V_p ^{\varepsilon}C}:\left(-\frac{1}{2q},\frac{1}{2q}\right)\times \V_p ^{\delta_p}(\IS^1p\subset C)\times \V_p ^{\varepsilon}C\rightarrow \V_p ^{\varepsilon}C
\end{align*}
with the restricted map $\Theta^{-1} _p$ has invertible differential at $0$.
Choose a sufficiently small neighborhood $U_p\subset \mathcal{U}_p$ of $p$ and a disc $D_p\subset \V_p ^{\varepsilon}C$ such that $P_{\V_p ^{\varepsilon}C}\circ \Theta_p ^{-1}$ restricts to a diffeomorphism from $U_p\cap T^1 _mM$ onto $D$. Denote this composition $\pr_{\V_p ^{\varepsilon}C}\circ \Theta^{-1} _p$ by
$$\theta_p:U_p \cap T^1 _mM\rightarrow D_p.$$
Then $\theta_p$ is a chart on $U_p\cap T^1 _mM$, and there are smooth functions
\begin{align*}
t_p:D_p\rightarrow \left(\frac{-1}{2q},\frac{1}{2q}\right)
\end{align*}
and
\begin{align*}
f_p:D_p\rightarrow \V_p ^{\delta_p}(\IS^1p\subset C),
\end{align*}
such that, in the charts $\Theta^{-1} _p$ and $\theta_p$, the inclusion $\iota:T^1 _mM\rightarrow T^1M$ is the graph map
\begin{align*}
\Theta^{-1} _p\circ \iota \circ \theta^{-1} _p:D_p&\rightarrow \Theta^{-1} _p(U_p\cap T^1 _mM)\subset \left(\frac{-1}{2q},\frac{1}{2q}\right)\times \V_p {(\IS^1p\subset C)}\times \V_p ^{\varepsilon}C,\\
x&\mapsto (t_p(x),f_p(x),x).
\end{align*}
By (\ref{transversalintersectionlagrangian}), $f_p$ is a diffeomorphism in a neighborhood of $\theta_p(U_p\cap T^1 _mM)$ after possibly shrinking $U_p$. In the charts $\theta_p$ and $\tilde{\phi}_p$, the map $\Pi_{\IS^1}$ is given by
\begin{align*}
\tilde{\phi}_p ^{-1}\circ \Pi_{\IS^1} \circ \theta^{-1} _p:D_p&\rightarrow \V_p ^{\delta_p}(\IS^1p\subset C)\times (B_{\varepsilon}(0)/(\IZ/q\IZ))\\
x&\mapsto (f_p(x),[ x]).
\end{align*}
Thus, in the charts $\theta_p$ and $\bar{\phi}_p$, the map $\Pi_{\IS^1}$ is given by
\begin{align*}
\bar{\phi}_p\circ \Pi_{\IS^1}\circ \theta_p ^{-1}:D_p &\rightarrow \V_p ^{\delta_p}(\IS^1p\subset C)\times B_{\varepsilon^q}(0),\\
x&\mapsto (f_p(x),x^{q}),
\end{align*}
which is a smooth embedding.
\end{proof}

\begin{lemma}\label{lemcharaccongruence}
Let $\V{\Pi_{\IS^1}}$ be the normal bundle of the immersion $\Pi_{\IS^1}|_{T_{m} ^1M}$. The euler characteristic of the normal bundle $\V{\Pi_{\IS^1}}$ satisfies $\chi(\V{\Pi_{\IS^1}})\equiv 2 \mod 4$.
\end{lemma}

\begin{proof}
On the preimage of the regular part of $T^1M/\IS^1$, an $\omega_{orb}$-compatible almost complex structure induces an isomorphism between $TT^1_mM$ and the normal bundle of $\Pi_{\IS^1}$. To prove Lemma \ref{lemcharaccongruence}, we will modify the normal bundle near the preimage of the singular set so that this isomorphism extends over $T^1_mM$. We will then show that the modification changes the normal Euler number by a multiple of four.\\

Let $C$ be a component of the isotropy set with isotropy group $\IZ/q\IZ$, and let $p\in T^1_mM\cap C$. By Proposition \ref{propsmoothimmersion}, there are a neighborhood $U_p$ of $p$, a disc $D_p\subset\V_p^{\varepsilon}C$ centered at $0$, a chart $\theta_p:T^1_mM\cap U_p\rightarrow D_p$, and a diffeomorphism $f_p:D_p\rightarrow\V_p^{\delta_p}(\IS^1p\subset C)$ onto its image such that
\begin{align*}
\tilde{\phi}_p^{-1}\circ\Pi_{\IS^1}\circ\theta_p^{-1}(x)&=(f_p(x),[x]),\\
\bar{\phi}_p^{-1}\circ\Pi_{\IS^1}\circ\theta_p^{-1}(x)&=(f_p(x),x^q).
\end{align*}
Recall that the orbifold chart is induced by
\begin{align*}
\phi_p:\V_p^{\delta_p}(\IS^1p\subset C)\times\V_p^{\varepsilon}C&\rightarrow T^1M,\\
(v,w)&\mapsto\exp(\mathcal{P}_{t\mapsto\exp(tv)}(w)).
\end{align*}
For sufficiently small $\delta_p$ and $\varepsilon$, our choice of metric ensures that $\tilde{\phi}_p^*\omega_{orb}$ lifts to $\phi_p^*d\bar{\alpha}$, which is the constant form
\[
((v_1,w_1),(v_2,w_2))\mapsto
\langle iv_1,v_2\rangle_{\IR^2}+\langle iw_1,w_2\rangle_{\IR^2}.
\]
The manifold coordinates are obtained from the orbifold coordinates by $(z_1,z_2)\mapsto(z_1,z_2^q)$. Away from $z_2=0$, this map is holomorphic and locally invertible. Thus multiplication by $i$ defines a $\bar{\phi}_p^*\omega_{orb}$-compatible almost complex structure there. We may therefore choose an $\omega_{orb}$-compatible almost complex structure $J$ on the regular part of $T^1M/\IS^1$ that is multiplication by $i$ in these charts. Since $\Pi_{\IS^1}$ is Lagrangian on this part, $J$ induces an isomorphism
\[
TT^1_mM|_{T^1_mM\setminus\bigcup_{i=1}^kC_i}
\rightarrow\V\Pi_{\IS^1}|_{T^1_mM\setminus\bigcup_{i=1}^kC_i}.
\]
We will use this local description to extend the induced bundle map after modifying $\V\Pi_{\IS^1}$. Also, the lift $x\mapsto(f_p(x),x)$ has Lagrangian image $Graph(f_p)\subset\IC^2$.\\

For this purpose, we regard the normal bundle as the quotient
\begin{align*}
\V\Pi_{\IS^1}
=\Pi_{\IS^1}^*T(T^1M/\IS^1)/d\Pi_{\IS^1}(TT^1_mM).
\end{align*}
Here $TT^1_mM$ is identified with its image in the pullback bundle under $v\mapsto(x,d\Pi_{\IS^1}(v))$, for $v\in T_x(T^1_mM)$. The charts $\theta_p$ and $\bar{\phi}_p$ identify $\Pi_{\IS^1}^*T(T^1M/\IS^1)|_{\theta_p^{-1}(D_p)}$ with $D_p\times\IC^2$ by
\begin{align*}
(x,v)\mapsto\bigl(\theta_p(x),(d\bar{\phi}_p^{-1})_{\Pi_{\IS^1}(x)}(v)\bigr).
\end{align*}
In these coordinates, the image of the tangent bundle is
\begin{align*}
TD_p=\{(z,(df_p(z)w,qz^{q-1}w)):z\in D_p,\ w\in\IC\},
\end{align*}
and the normal bundle is
\begin{align*}
\V D_p=(D_p\times\IC^2)/TD_p.
\end{align*}
We now choose trivializations of these bundles and compute the map induced by $J$. The assignment
\begin{align*}
\psi_p:(z,(df_p(z)w,qz^{q-1}w))\mapsto(z,w)
\end{align*}
trivializes $TD_p$. Each class in
\[
\V_zD_p=\IC^2/\{(df_p(z)w,qz^{q-1}w):w\in\IC\}
\]
has a unique representative with first component zero. For $[v,w]$, this representative is $(0,w-qz^{q-1}df_p(z)^{-1}v)$, so
\begin{align*}
\Psi_p:\V D_p&\rightarrow D_p\times\IC,\\
(z,[v,w])&\mapsto(z,w-qz^{q-1}df_p(z)^{-1}v)
\end{align*}
is a trivialization. Consequently, the bundle map induced by $J$ is given by
\begin{align}\label{complexstructureintrivialisations}
\Psi_p\circ J\circ\psi_p^{-1}:D_p\times\IC&\rightarrow D_p\times\IC\nonumber\\
(z,w)&\mapsto(z,qz^{q-1}(iw-df_p(z)^{-1}idf_p(z)w)).
\end{align}

We now show that $w\mapsto i\cdot w-df_p(z)^{-1}i \cdot df_p(z)w$ is an isomorphism for every $z\in D_p$. By Proposition \ref{propsmoothimmersion}, $Graph(f_p)$ is Lagrangian for the standard symplectic form displayed above. Its tangent bundle is
\begin{align*}
TGraph(f_p)=\{((f_p(z),z),(df_p(z)w,w)):z\in D_p,\ w\in\IC\},
\end{align*}
with trivialization
\begin{align*}
\bar{\psi}_p:TGraph(f_p)&\rightarrow Graph(f_p)\times\IC,\\
((f_p(z),z),(df_p(z)w,w))&\mapsto((f_p(z),z),w).
\end{align*}
The normal bundle is the quotient of $Graph(f_p)\times\IC^2$ by $TGraph(f_p)$. Since $f_p$ is a diffeomorphism onto its image, each normal class $[v,w]$ has the unique representative $(0,w-df_p(z)^{-1}v)$ with first component zero. Thus
\begin{align*}
\bar{\Psi}_p:\V Graph(f_p)&\rightarrow Graph(f_p)\times\IC,\\
((f_p(z),z),[v,w])&\mapsto((f_p(z),z),w-df_p(z)^{-1}v)
\end{align*}
is a trivialization. Multiplication by $i$ maps the tangent space of the Lagrangian graph to its orthogonal complement and hence induces a bundle isomorphism
\[
J:TGraph(f_p)\rightarrow\V Graph(f_p).
\]
In these trivializations, it is given by
\begin{align}\label{liftiso}
\bar{\Psi}_p\circ J\circ\bar{\psi}_p^{-1}:Graph(f_p)\times\IC&\rightarrow Graph(f_p)\times\IC\nonumber\\
((f_p(z),z),w)&\mapsto((f_p(z),z),iw-df_p(z)^{-1}idf_p(z)w).
\end{align}
This proves the claimed invertibility.\\

By (\ref{complexstructureintrivialisations}), the map induced by $J$ is an isomorphism for $z\neq0$, but the factor $z^{q-1}$ prevents it from remaining an isomorphism at $z=0$. We remove this factor by modifying the normal bundle near $p$. Identify $D_p$ with $\theta_p^{-1}(D_p)$ when discussing bundle restrictions, and set
\[
\frac12D_p\coloneqq\{z/2:z\in D_p\}.
\]
We cut out the bundle over the closed half-disc $\overline{\frac12D_p}$ and glue in $D_p\times\IC$ using
\[
(z,v)\sim\Psi_p^{-1}(z,z^{q-1}v),
\qquad z\in D_p\setminus\overline{\frac12D_p}.
\]
This gives the modified bundle $(\V\Pi_{\IS^1})_{mod,p}$. A trivialization $\Phi_p$ over $D_p$ is given on the overlap and on the inserted bundle, respectively, by
\begin{align*}
\Phi_p([\Psi_p^{-1}(z,v)])&=(z,z^{-(q-1)}v),\\
\Phi_p([(z,v)])&=(z,v),
\end{align*}
where brackets denote classes under the gluing identification. On the overlap, the induced map is therefore
\begin{align*}
\Phi_p\circ J\circ\psi_p^{-1}:\left(D_p\setminus\overline{\frac12D_p}\right)\times\IC
&\rightarrow\left(D_p\setminus\overline{\frac12D_p}\right)\times\IC,\\
(z,w)&\mapsto(z,z^{-(q-1)}qz^{q-1}(iw-df_p(z)^{-1}idf_p(z)w))\\
&=(z,q(iw-df_p(z)^{-1}idf_p(z)w)).
\end{align*}
By (\ref{liftiso}), this extends to a bundle isomorphism over all of $D_p$. Repeating the construction in disjoint neighborhoods of the points of $T^1_mM\cap\bigcup_{i=1}^kC_i$ gives a bundle $(\V\Pi_{\IS^1})_{mod}$ isomorphic to $TT^1_mM$.\\

It remains to show that the modifications change the Euler number by a multiple of four. We write
\[
T^1_mM\cap\bigcup_{i=1}^kC_i=\{p_1,\ldots,p_r\},
\qquad p_j\in C_{i_j},
\]
and set $(\V\Pi_{\IS^1})^0=\V\Pi_{\IS^1}$. For $j=1,\ldots,r$, let $(\V\Pi_{\IS^1})^j$ be the bundle obtained from $(\V\Pi_{\IS^1})^{j-1}$ by the modification at $p_j$. Thus $(\V\Pi_{\IS^1})^r=(\V\Pi_{\IS^1})_{mod}$. Since the modification neighborhoods are disjoint, the trivialization $\Psi_{p_j}$ remains available over $D_{p_j}$ before the $j$-th modification.

We now fix $j$. Since $T^1_mM\cong\IS^2$, the complement of the closed half-disc at $p_j$ is an open disc, over which $(\V\Pi_{\IS^1})^{j-1}$ is trivial. We choose a trivialization $\Psi^{j-1}$ there. Its domain overlaps with that of the local trivialization $\Psi_{p_j}$ in the annulus
\[
D_{p_j}\setminus\overline{\frac12D_{p_j}}.
\]
After choosing orientations, the transition map takes values in the orientation-preserving linear automorphisms of $\IR^2$. Restricting this map to a circle in the annulus and deforming it to $SO(2)\cong\IS^1$ gives a map
\[
G^{j-1}_{p_j}:\IS^1\rightarrow\IS^1.
\]
With consistent orientation conventions, the Euler number is given by
\begin{align*}
\chi((\V\Pi_{\IS^1})^{j-1})=\deg G^{j-1}_{p_j}.
\end{align*}

The modification leaves the bundle unchanged outside the half-disc, so we retain the exterior trivialization $\Psi^{j-1}$. On the overlap, the new local trivialization $\Phi_{p_j}$ satisfies
\[
\Phi_{p_j}\bigl([\Psi_{p_j}^{-1}(z,v)]\bigr)
=(z,z^{-(q_{i_j}-1)}v).
\]
Thus the transition map changes by multiplication by $z^{-(q_{i_j}-1)}$ in these local coordinates. On a circle $z=\rho e^{it}$ in the annulus, this factor takes the form
\[
z^{-(q_{i_j}-1)}
=\rho^{-(q_{i_j}-1)}e^{-i(q_{i_j}-1)t}.
\]
The positive real factor does not affect the degree, whereas the rotation has winding number $-(q_{i_j}-1)$. Taking into account the orientations of the local base and fiber coordinates, the transition map $G^j_{p_j}$ of the modified bundle therefore satisfies
\begin{align*}
\chi((\V\Pi_{\IS^1})^j)
&=\deg G^j_{p_j}=\deg G^{j-1}_{p_j}\pm(q_{i_j}-1)\\
&=\chi((\V\Pi_{\IS^1})^{j-1})\pm(q_{i_j}-1).
\end{align*}
Here the bundle orientations are chosen to agree outside the modification region. Summing the changes over all modification points gives
\begin{align*}
\chi((\V\Pi_{\IS^1})_{mod})-\chi(\V\Pi_{\IS^1})
=\sum_{i=1}^k\sum_{p\in T^1_mM\cap C_i}\pm(q_i-1).
\end{align*}

For each $i$, the finite set $T^1_mM\cap C_i$ is invariant under $v\mapsto-v$. Since this involution has no fixed points on $T^1_mM$, the points occur in pairs, and the set therefore has even cardinality. Moreover, $q_i$ is odd, so $q_i-1$ is even. This implies that
\[
\pm(q_i-1)\equiv q_i-1\mod4,
\]
so the signs of the individual contributions do not affect the sum modulo four. For each component $C_i$, we therefore obtain
\begin{align*}
\sum_{p\in T^1_mM\cap C_i}\pm(q_i-1)
&\equiv \#(T^1_mM\cap C_i)(q_i-1)\\
&\equiv0\mod4.
\end{align*}
Summing over the components, we conclude that
\begin{align*}
\chi((\V\Pi_{\IS^1})_{mod})\equiv\chi(\V\Pi_{\IS^1})\mod4.
\end{align*}
Finally, $(\V\Pi_{\IS^1})_{mod}\cong TT^1_mM$, whose Euler number is $2$ with the tangent orientation. The bundle isomorphism may reverse the chosen fiber orientation, so the Euler number of the modified bundle is $\pm2$. In either case, we obtain
\begin{align*}
\chi(\V\Pi_{\IS^1})\equiv\chi((\V\Pi_{\IS^1})_{mod})\equiv2\mod4.
\end{align*}
This completes the proof.
\end{proof}

For an embedding, the normal Euler number equals the homological self-intersection, so Lemma \ref{lemcharaccongruence} would already give the desired class. For an immersion, we instead apply the self-intersection formula to a suitable $C^2$-approximation with only transverse double points.


\subsection{Topological consequences for the quotient \texorpdfstring{$T^1M/\IS^1$}{T1MS1}}

We can now determine the homeomorphism type of the quotient $T^1M/\IS^1$.

\begin{proof}[Proof of Theorem \ref{quotclasstheorem}]
By Proposition \ref{propquotienthomeomorphismtyp}, it remains to find a class $c\in H_2(T^1M/\IS^1;\IZ)$ with self-intersection $c^2\equiv 2\mod 4$. We combine the classical self-intersection formula with the results of the preceding subsection. Choose $m\in M$ such that $T^1_mM$ intersects each $C\in\{C_1,\ldots,C_k\}$ transversally, and equip $T^1M/\IS^1$ with the smooth atlas $\hat{\mathcal{A}}$. By Proposition \ref{propsmoothimmersion}, the restriction $\Pi_{\IS^1}:T^1_mM\rightarrow T^1M/\IS^1$ is a smooth immersion, and Lemma \ref{lemcharaccongruence} gives $\chi(\V\Pi_{\IS^1})\equiv 2\mod 4$. For a generator $[T^1_mM]\in H_2(T^1_mM;\IZ)$, we will show that
\begin{align*}
(\Pi_{\IS^1})_*([T^1_mM])\cdot(\Pi_{\IS^1})_*([T^1_mM])\equiv 2\mod 4.
\end{align*}
For this, we use the antipodal symmetry and the density, in the $C^2$-topology, of immersions with only finitely many transverse double points (see \cite{wall2016differential}) to construct a suitable approximation $f$ of $\Pi_{\IS^1}$.\\

The inclusion $T^1_mM\rightarrow T^1M$ is equivariant with respect to the antipodal map $A$. Thus $\Pi_{\IS^1}$ is equivariant with respect to $A$ and the induced smooth fixed-point-free involution $\bar{A}$ on $T^1M/\IS^1$. The quotient $(T^1M/\IS^1)/(\IZ/2\IZ)$ is therefore a smooth manifold, and the projection
\[
\Pi_{\IZ/2\IZ}:T^1M/\IS^1\rightarrow (T^1M/\IS^1)/(\IZ/2\IZ)
\]
is a local diffeomorphism. We use the same notation for the antipodal quotient projection of $T^1_mM$. Equivariance induces a smooth map
\begin{align*}
\bar{\Pi}_{\IS^1}:T^1_mM/(\IZ/2\IZ)\rightarrow (T^1M/\IS^1)/(\IZ/2\IZ)
\end{align*}
such that the following diagram commutes.
\begin{center}
\begin{tikzpicture}[commutative diagrams/every diagram]
\matrix[matrix of math nodes, name=m, column sep= 1.5cm, commutative diagrams/every cell] {
T^1 _mM  & T^1M/\IS^1 \\
T^1 _m M/(\IZ/2\IZ)& (T^1M/\IS^1)/(\IZ/2\IZ) \\};
\path[commutative diagrams/.cd, every arrow, every label]
(m-1-1) edge node{$\Pi_{\IS^1}$} (m-1-2)
(m-2-1)edge node{$\bar{\Pi}_{\IS^1}$}(m-2-2)
(m-1-2) edge node{$\Pi_{\IZ/2\IZ}$} (m-2-2)
(m-1-1) edge node{$\Pi_{\IZ/2\IZ}$} (m-2-1);
\end{tikzpicture}
\end{center}
Since both vertical maps are local diffeomorphisms, $\bar{\Pi}_{\IS^1}$ is an immersion. Choose a $\IZ/2\IZ$-invariant metric on $T^1M/\IS^1$ and equip its quotient with the metric for which $\Pi_{\IZ/2\IZ}$ is a local isometry. Its differential
\begin{align*}
d(\Pi_{\IZ/2\IZ})_p:T_p(T^1M/\IS^1)\rightarrow T_{\Pi_{\IZ/2\IZ}(p)}\bigl((T^1M/\IS^1)/(\IZ/2\IZ)\bigr)
\end{align*}
is then a linear isometry. In particular, the map
\begin{align*}
Q:\Pi_{\IS^1}^*T(T^1M/\IS^1)&\rightarrow \bar{\Pi}_{\IS^1}^*T\bigl((T^1M/\IS^1)/(\IZ/2\IZ)\bigr),\\
(x,v)&\mapsto\bigl(\Pi_{\IZ/2\IZ}(x),d(\Pi_{\IZ/2\IZ})_{\Pi_{\IS^1}(x)}(v)\bigr)
\end{align*}
is a bundle map covering the antipodal quotient projection of $T^1_mM$ and is an isometry on each fiber.\\

If a smooth map $\bar{f}$ is sufficiently close to $\bar{\Pi}_{\IS^1}$ in the $C^2$-topology, it is an immersion and can be written as
\[
\bar{f}(\bar{x})=\exp_{\bar{\Pi}_{\IS^1}(\bar{x})}\bigl(\bar{\eta}(\bar{x})\bigr)
\]
for a vector field $\bar{\eta}$ along $\bar{\Pi}_{\IS^1}$ with small $C^2$-norm. The fiberwise isometry $Q$ determines a unique vector field $\eta$ along $\Pi_{\IS^1}$ satisfying
\[
d(\Pi_{\IZ/2\IZ})_{\Pi_{\IS^1}(x)}\bigl(\eta(x)\bigr)
=\bar{\eta}\bigl(\Pi_{\IZ/2\IZ}(x)\bigr).
\]
We define
\[
f(x)=\exp_{\Pi_{\IS^1}(x)}\bigl(\eta(x)\bigr).
\]
Since $\Pi_{\IZ/2\IZ}$ is a Riemannian covering, $f$ satisfies $\Pi_{\IZ/2\IZ}\circ f=\bar{f}\circ\Pi_{\IZ/2\IZ}$ and is equivariant with respect to $A$ and $\bar{A}$.

To compare the normal bundles, we use parallel translation along the curves
\[
t\mapsto\exp_{\Pi_{\IS^1}(x)}\bigl((1-t)\eta(x)\bigr),
\qquad
t\mapsto\exp_{\bar{\Pi}_{\IS^1}(\bar{x})}\bigl((1-t)\bar{\eta}(\bar{x})\bigr).
\]
This gives bundle isomorphisms
\begin{align*}
f^*T(T^1M/\IS^1)&\rightarrow\Pi_{\IS^1}^*T(T^1M/\IS^1),\\
\bar{f}^*T\bigl((T^1M/\IS^1)/(\IZ/2\IZ)\bigr)
&\rightarrow\bar{\Pi}_{\IS^1}^*T\bigl((T^1M/\IS^1)/(\IZ/2\IZ)\bigr).
\end{align*}
They identify $\V f$ and $\V\bar{f}$ with subbundles of the corresponding pullback bundles. If the $C^2$-norm of $\bar{\eta}$ is sufficiently small, parallel translation followed by orthogonal projection onto $\V\bar{\Pi}_{\IS^1}$ gives an isomorphism $\V\bar{f}\cong\V\bar{\Pi}_{\IS^1}$. Since the covering preserves the metric and parallel translation, the same construction gives $\V f\cong\V\Pi_{\IS^1}$.\\

We can therefore choose $\bar{f}$ with only finitely many transverse double points and sufficiently close to $\bar{\Pi}_{\IS^1}$ that these normal-bundle isomorphisms hold. Its equivariant lift $f$ also has only finitely many transverse double points. Indeed, equivariance and the absence of fixed points of $\bar{A}$ imply that $f(x)\neq f(-x)$, while both quotient projections are local diffeomorphisms. If $f(x)=f(y)=p$, then
\[
f(-x)=f(-y)=\bar{A}(p).
\]
Thus the double points occur in distinct pairs $p,\bar{A}(p)$. Each pair contributes either $\pm2$ or $0$ to the signed double-point count $I(f)$, so $I(f)$ is even.\\


By the self-intersection formula for immersed oriented surfaces \cite[p.~215]{yang1984self}, we obtain
\begin{align*}
f_*([T^1_mM])\cdot f_*([T^1_mM])
&=\chi(\V f)+2I(f)\\
&\equiv 2\mod 4,
\end{align*}
where each double point is counted once and $\chi(\V f)$ is the normal Euler number. The congruence follows from the evenness of $I(f)$, the isomorphism $\V f\cong\V\Pi_{\IS^1}$, and Lemma \ref{lemcharaccongruence}. Since the perturbation $f$ is homotopic to $\Pi_{\IS^1}$, this proves the claim.

\end{proof}

%% file: maintheoremproof.tex
\section{Proof of the main theorem}

\begin{lemma}\label{lembasisexistence}
Let $(M,\langle \cdot,\cdot \rangle)$ be a simply connected Riemannian manifold all of whose geodesics are closed. Suppose that the $\IS^1$-action on the unit tangent bundle induced by the geodesic flow has a nonempty isotropy set. Then one can choose a basis $E_1,E_2$ of $H_2(T^1M/\IS^1;\IZ)$ and, for $i=1,\ldots,k$, integers $l_i\in\IZ$ such that
\[
E_1^2=E_2^2=0,\qquad E_1\cdot E_2=1,
\]
and
\[
[C_i/\IS^1]=E_1+l_iE_2,\qquad i=1,\ldots,k.
\]
If $k>2$, then $l_i=0$ for every $i=1,\ldots,k$.
\end{lemma}

\begin{proof}
The non-free locus of the $\IS^1$-action induced by the geodesic flow consists of finitely many pairwise disjoint connected $3$-manifolds $C_1,\ldots,C_k$. By Theorem B of the authors prior paper "Isotropy sets of periodic geodesic flows on three-dimensional manifolds", there are at least two such components, and each $C_i$ is the fixed-point set of a finite cyclic subgroup $\IZ/q_i\IZ$. By Theorem \ref{propquotienthomeomorphismtyp}, the quotient $T^1M/\IS^1$ carries a smooth structure $\mathcal{B}$ and a smooth symplectic form $\omega$ that agree with the corresponding orbifold structures away from $\bigcup_i C_i/\IS^1$. Furthermore, each $C_i/\IS^1$ is a symplectic surface, with
\[
\omega|_{C_i/\IS^1}=\omega_{orb}|_{C_i/\IS^1}
\]
as area forms. By Theorem \ref{propquotienthomeomorphismtyp}, $T^1M/\IS^1$ is homeomorphic to $\IS^2\times\IS^2$ and is therefore minimal. Since the surfaces $C_i/\IS^1$ are pairwise disjoint, their homology classes are pairwise orthogonal. It follows that at least one of the classes $[C_1/\IS^1]$ and $[C_2/\IS^1]$ has nonnegative self-intersection. Relabeling if necessary, we may assume that
\[
[C_1/\IS^1]^2\geq 0.
\]
McDuff's theorem (Theorem 1.4 in \cite{Mcduff1990TheSO}) then provides a compact Riemann surface $\Sigma$, a symplectic sphere bundle $\pi_P:P\to\Sigma$, and a diffeomorphism
\[
F:P\longrightarrow (T^1M/\IS^1,\mathcal{B})
\]
such that $F^{-1}(C_1/\IS^1)$ is a symplectic section. As $P$ is homeomorphic to $\IS^2\times\IS^2$, the long exact sequence of the fibration shows that $\Sigma$ is simply connected, and thus diffeomorphic to $\IS^2$. Among the two smooth $\IS^2$-bundles over $\IS^2$, only the trivial one has total space homeomorphic to $\IS^2\times\IS^2$. Consequently, we may regard
\[
F:\IS^2\times\IS^2\longrightarrow (T^1M/\IS^1,\mathcal{B})
\]
as a diffeomorphism for which the fibers $\{p\}\times\IS^2$ are symplectic with respect to $F^*\omega$. If $k>2$, then $[C_i/\IS^1]^2=0$ for every $i=1,\ldots,k$. Interchanging $E_1$ and $E_2$ if necessary, we may therefore take $[C_1/\IS^1]=E_1$. Pairwise orthogonality then yields
\[
[C_i/\IS^1]=E_1,\qquad i=1,\ldots,k.
\]
Hence $l_i=0$ for all $i$. In the case $k=2$, there is a positive integer $l$ such that
\[
[C_1/\IS^1]=E_1+lE_2,\qquad [C_2/\IS^1]=E_1-lE_2.
\]
\end{proof}

\begin{proof}[Proof of Theorem~\ref{bergconj}]
Suppose, for a contradiction, that the geodesic-flow action has a nonempty nontrivial isotropy set, with connected components
\[
C_1,\ldots,C_k.
\]
Theorem B of the authors prior paper "Isotropy sets of periodic geodesic flows on three-dimensional manifolds" then immediately implies that $k\geq 2$. We complete the argument outlined at the beginning of Chapter 3. By Satake \cite{satake1956generalization}, smooth singular homology agrees with ordinary singular homology over $\IR$, and orbifold de Rham cohomology agrees with singular cohomology. Under these identifications, the natural pairing is given by integration:
\begin{align}\label{eqderhamisoorb}
\langle[\omega],[c]\rangle=\int_c\omega.
\end{align}
In particular, $\omega_{orb}$ determines a class $[\omega_{orb}]\in H^2(T^1M/\IS^1;\IR)$. By Weinstein's argument (Proposition B in \cite{weinstein1977symplectic}), $lcm(q_1,\ldots,q_k)[\omega_{orb}]$ is an integral class.

Since the action has more than one orbit type, Lemma \ref{lembasisexistence} gives, for $i=1,...,k$, integers $l_i\in \IZ$ and a basis $E_1,E_2$ of $H_2(T^1M/\IS^1;\IZ)$ satisfying
\[
E_1^2=E_2^2=0,\qquad E_1\cdot E_2=1,\qquad [C_i/\IS^1]=E_1+l_iE_2.
\]
Let $E_1^*,E_2^*$ denote the dual basis of $H^2(T^1M/\IS^1;\IZ)$; we retain this notation after extension of scalars to $\IR$. Write
\[
[\omega_{orb}]=\omega_1E_1^*+\omega_2E_2^*.
\]
Consequently, $lcm(q_1,\ldots,q_k)\omega_1,lcm(q_1,\ldots,q_k)\omega_2\in\IZ$. We shall first show that both integers are odd and then derive a contradiction by computing $\omega_1$ from the surface $C_1/\IS^1$.\\

\noindent Choose $m\in M$ such that $T^1_mM$ meets every $C_i$ transversely, and set
\begin{align}\label{basisrepresentation}
A\coloneqq(\Pi_{\IS^1})_*([T^1_mM])=aE_1+bE_2.
\end{align}
The intersection form is a topological invariant. For $U,V\in H_2(T^1M/\IS^1;\IZ)$, it is given by
\[
U\cdot V=\langle PD^{-1}(U)\cup PD^{-1}(V),[T^1M/\IS^1]\rangle.
\]
In the Proof from Theorem B it was shown that we have the congruence
\[
A^2\equiv 2\pmod 4.
\]
On the other hand, equation \eqref{basisrepresentation} gives $A^2=2ab$, so
\[
a\equiv b\equiv 1\pmod 2.
\]
Moreover, Proposition \ref{propsmoothimmersion} shows that the map $\Pi_{\IS^1}:T^1_mM\to T^1M/\IS^1$ is a Lagrangian orbifold immersion. Hence \eqref{eqderhamisoorb} yields
\[
0=\int_{T^1_mM}\Pi_{\IS^1}^*\omega_{orb}
=[\omega_{orb}](A)
=\omega_1a+\omega_2b.
\]
In particular,
\begin{align}\label{equationvanishingsymplecticform}
\omega_1a+\omega_2b=0.
\end{align}
Multiplying \eqref{equationvanishingsymplecticform} by $lcm(q_1,\ldots,q_k)$ and reducing modulo $2$, we obtain
\[
lcm(q_1,\ldots,q_k)\omega_1\equiv lcm(q_1,\ldots,q_k)\omega_2\pmod 2,
\]
because $a$ and $b$ are odd. By Theorem 5.5 of the authors prior paper "Isotropy sets of periodic geodesic flows on three-dimensional manifolds", each $q_i$ is odd, so $lcm(q_1,\ldots,q_k)$ is odd as well.
We claim that
\[
lcm(q_1,\ldots,q_k)\omega_1\equiv lcm(q_1,\ldots,q_k)\omega_2\equiv 1\pmod 2.
\]
Suppose instead that both integers are even, and write
\[
lcm(q_1,\ldots,q_k)\omega_1=2b_1,\qquad lcm(q_1,\ldots,q_k)\omega_2=2b_2
\]
for $b_1,b_2\in\IZ$. Let
\[
Q:T^1M\times_{\IS^1} E\IS^1\longrightarrow T^1M/\IS^1
\]
be the natural projection. The groups $H_{\IS^1}^*(T^1M;\IZ)$ are finitely generated and are free in degrees at most $3$. Choose an integral basis $F_1^*,F_2^*$ of $H_{\IS^1}^2(T^1M;\IZ)$. Since $b_1E_1^*+b_2E_2^*$ is an integral class, we may write
\[
Q^*(b_1E_1^*+b_2E_2^*)=B_1F_1^*+B_2F_2^*
\]
with $B_1,B_2\in\IZ$. As $Q$ is a rational homotopy equivalence, it induces an isomorphism
\[
Q^*:H^*(T^1M/\IS^1;\IR)\longrightarrow H_{\IS^1}^*(T^1M;\IR).
\]
By orbifold Chern--Weil theory (see p.~120 of \cite{boyer2007sasakian}), the image $e_{\IR}$ of the integral Euler class $e_{\IZ}\in H_{\IS^1}^2(T^1M;\IZ)$ in real cohomology satisfies
\[
e_{\IR}=Q^*[\omega_{orb}].
\]
Substituting $\omega_i=2b_i/lcm(q_1,\ldots,q_k)$, we obtain
\[
e_{\IR}
=\frac{2}{lcm(q_1,\ldots,q_k)}Q^*(b_1E_1^*+b_2E_2^*)
=\frac{2B_1}{lcm(q_1,\ldots,q_k)}F_1^*+\frac{2B_2}{lcm(q_1,\ldots,q_k)}F_2^*.
\]
Since $e_{\IZ}$ is integral, $2B_i/lcm(q_1,\ldots,q_k)\in\IZ$. As $lcm(q_1,\ldots,q_k)$ is odd, this implies $B_i/lcm(q_1,\ldots,q_k)\in\IZ$, and hence
\[
e_{\IZ}
=2\left(\frac{B_1}{lcm(q_1,\ldots,q_k)}F_1^*+\frac{B_2}{lcm(q_1,\ldots,q_k)}F_2^*\right).
\]
Thus $e_{\IZ}$ is divisible by two and therefore not primitive. This contradicts the Gysin sequence of the principal $\IS^1$-bundle
\[
T^1M\times E\IS^1\longrightarrow T^1M\times_{\IS^1} E\IS^1,
\]
because
\[
H^2(T^1M\times E\IS^1;\IZ)\cong H^2(T^1M;\IZ)\cong\IZ
\]
is torsion-free. Hence
\[
lcm(q_1,\ldots,q_k)\omega_1\equiv lcm(q_1,\ldots,q_k)\omega_2\equiv 1\pmod 2.
\]
We distinguish the cases $k=2$ and $k>2$, starting with $k>2$. In this case, $[C_i/\IS^1]=E_1$ for $i=1,...,k$. We compute $\omega_1$ by evaluating on $C_1/\IS^1$. The class
\[
q_1[\omega_{orb}|_{C_1/\IS^1}]
\]
is the image of the integral Euler class $e(C_1\to C_1/\IS^1)\in H^2(C_1/\IS^1;\IZ)$ of the principal $\IS^1$-bundle
\[
\Pi_{\IS^1}:C_1\longrightarrow C_1/\IS^1
=(C_1/(\IZ/q_1\IZ))/(\IS^1/(\IZ/q_1\IZ)).
\]
Let $B\in H_2(C_1/\IS^1;\IZ)$ be the generator corresponding to the orientation induced by $\omega_{orb}|_{C_1/\IS^1}$. Since $C_1/\IS^1$ is homeomorphic to $\IS^2$, the Gysin sequence gives
\[
q_1[\omega_{orb}](E_1)
=e(C_1\to C_1/\IS^1)(B)
=\#\pi_1(C_1).
\]
Consequently,
\[
\omega_1=[\omega_{orb}](E_1)=\frac{\#\pi_1(C_1)}{q_1},
\]
and hence
\[
q_1\,lcm(q_1,\ldots,q_k)\omega_1=lcm(q_1,\ldots,q_k)\cdot\#\pi_1(C_1).
\]
The left-hand side is odd, since both $q_1$ and $lcm(q_1,\ldots,q_k)\omega_1$ are odd, whereas the right-hand side is even, since $\#\pi_1(C_1)$ is even. This contradiction shows that $k$ cannot exceed two.\\

\bigskip
\noindent It remains to rule out $k=2$. Let $L>0$ be the least common period of all unit-speed geodesics. Then $\phi_L=\mathrm{id}_{T^1M}$. Let $K_{i_0}$ denote the full critical manifold of the energy functional at energy $L^2/2$; under the initial-velocity identification, $K_{i_0}$ is diffeomorphic to $T^1M$.
We show that the orbifold first Chern class is a multiple of the Euler class and compute this multiple in two ways: first using a principal geodesic, and then using the exceptional components. We divide the argument into four claims.

\medskip
\noindent\textbf{Claim 1.}
The orbifold first Chern class is a multiple of the real Euler class. Under the real cohomology identification induced by $Q$, this means that for some $\kappa\in\IR$,
\begin{equation}\label{eqtwocomponentchernmultiplelcm}
c_1^{\mathrm{orb}}\bigl(T(T^1M/\IS^1)\bigr)
=\kappa[\omega_{orb}].
\end{equation}
Let $\bar\alpha$ be the Liouville form restricted to $T^1M$,
\[
\bar\alpha_{(x,v)}(Z)=\langle v,\D\pi(Z)\rangle,
\]
and put
\[
\xi=\ker\bar\alpha,
\qquad
\mathcal V=\ker\D\pi.
\]
With the chosen normalization of $\omega_{orb}$,
\[
\Pi_{\IS^1}^*\omega_{orb}=\frac1L\D\bar\alpha,
\qquad
\frac1L\int_\Gamma\bar\alpha=1
\]
for every principal orbit $\Gamma$. The vertical bundle $\mathcal V$ is an oriented Lagrangian subbundle of $(\xi,\D\bar\alpha)$. Its orientation is induced by the orientation of $M$, using $\mathcal V_{(x,v)}\cong v^\perp$.
Choose a flow-invariant compatible complex structure $\mathcal J$ on $\xi$, and write
\[
g_{\mathcal J}(u,v)=\D\bar\alpha(u,\mathcal Jv).
\]
Such an invariant choice exists because the geodesic flow defines a circle action. 
The Lagrangian property and compatibility of $\mathcal J$ give the real decomposition
\[
\xi=\mathcal V\oplus\mathcal J\mathcal V.
\]
Indeed, if $u\in\mathcal V$ and $\mathcal Ju\in\mathcal V$, then
\[
g_{\mathcal J}(u,u)=\D\bar\alpha(u,\mathcal Ju)=0,
\]
so $u=0$. This decomposition identifies $\xi$ with the complexification of $\mathcal V$ through the complex bundle isomorphism
\[
\mathcal V\otimes_{\IR}\IC\longrightarrow\xi,
\qquad
a+ib\longmapsto a+\mathcal Jb.
\]
The orientation and the metric $g_{\mathcal J}|_{\mathcal V}$ define rotation by $\pi/2$ in each fiber, making $\mathcal V$ a complex line bundle $\mathcal L$. Its complexification decomposes as
\[
\mathcal V\otimes_{\IR}\IC
\cong\mathcal L\oplus\overline{\mathcal L}.
\]
Here $\overline{\mathcal L}$ denotes the conjugate complex line bundle. After complexification, rotation by $\pi/2$ has eigenvalues $i$ and $-i$; the corresponding eigenbundles are isomorphic to $\mathcal L$ and $\overline{\mathcal L}$. Thus, as complex bundles,
\[
\xi\cong\mathcal L\oplus\overline{\mathcal L}.
\]
By additivity of the first Chern class and the identity $c_1(\overline{\mathcal L})=-c_1(\mathcal L)$, we obtain
\[
c_1(\xi)
=c_1(\mathcal L)+c_1(\overline{\mathcal L})
=c_1(\mathcal L)-c_1(\mathcal L)
=0.
\]
To define the orbifold first Chern class, consider the classifying space
\[
T^1M\times_{\IS^1} E\IS^1
\]
and the ordinary complex vector bundle
\[
\mathscr E=\xi\times_{\IS^1} E\IS^1
\longrightarrow T^1M\times_{\IS^1} E\IS^1.
\]
By definition,
\[
c_1^{\mathrm{orb}}\bigl(T(T^1M/\IS^1)\bigr)
:=c_1(\mathscr E).
\]
Fix $a_0\in E\IS^1$. The map
\[
j:T^1M\longrightarrow T^1M\times_{\IS^1} E\IS^1,
\qquad
j(u)=[u,a_0],
\]
satisfies $j^*\mathscr E\cong\xi$, and therefore $j^*c_1(\mathscr E)=0$.
The Gysin sequence of the principal circle bundle
\[
T^1M\times E\IS^1\longrightarrow T^1M\times_{\IS^1} E\IS^1
\]
shows that the kernel of this pullback in degree two is spanned by its Euler class. Thus $c_1(\mathscr E)$ is a multiple of the real Euler class $e_{\IR}=Q^*[\omega_{orb}]$, as identified in the first part of the proof.
Identifying the real cohomology of the classifying space with orbifold de Rham cohomology gives \eqref{eqtwocomponentchernmultiplelcm}.

\medskip
\noindent\textbf{Claim 2.}
For a principal geodesic, the determinant winding of the linearized flow, expressed in the disk frame below, is
\[
w=q_1+q_2.
\]
We first determine the index of $K_{i_0}$. As discussed in the author's earlier paper, \emph{Isotropy sets of periodic geodesic flows on three-dimensional manifolds} (in the section on index monotonicity), for each critical manifold $K_i$ consisting of exceptional geodesics, the next critical manifold $K_{i+1}$ at higher energy has index
\[
\ind(K_{i+1})=\ind(K_i)+2.
\]
Since the critical manifold of lowest positive energy has index $2$, the critical manifold $K_{i_0}$ at energy $L^2/2$ has index
\begin{equation}\label{eqtwocomponentindexlcm}
\ind(K_{i_0})=2(q_1+q_2)-2.
\end{equation}
Let $\gamma:[0,L]\to M$ be a unit-speed regular geodesic of energy $L^2$, and let
\[
\Gamma(t)=\bigl(\gamma(t),\dot\gamma(t)\bigr)
\]
be its lifted orbit. Since $T^1M$ is simply connected, there is a smooth map
\[
D:D_+\longrightarrow T^1M
\]
from an oriented disk whose positively oriented boundary parametrizes $\Gamma$.\\

\noindent The bundle $D^*\mathcal V$ is an oriented real vector bundle of rank two over a contractible disk, so it is trivial. Choose an oriented $g_{\mathcal J}$-orthonormal frame $e_1,e_2$ of $D^*\mathcal V$. It determines a unitary frame $e$ of $D^*\xi$ by
\[
e(a+ib)
=a_1e_1+a_2e_2+b_1\mathcal Je_1+b_2\mathcal Je_2.
\]
Thus $e(\IR^2)=\mathcal V$. Along the boundary define
\[
G(t)=e(t)^{-1}(\D\phi_t)_{\Gamma(0)}e(0).
\]
The flow preserves $\D\bar\alpha$ and $\mathcal J$, so $G(t)\in U(2)$. Since $\phi_L=\mathrm{id}$ and the frame is periodic,
\[
G(0)=G(L)=I.
\]
Set
\begin{equation}\label{eqtwocomponentdetwinding}
w=\deg(\det G).
\end{equation}
We relate this degree to Jacobi fields using the connection-induced splitting.
The Levi--Civita connection gives a splitting of the tangent bundle of the unit tangent bundle into horizontal and vertical subbundles,
\[
T(T^1M)=\mathcal H\oplus\mathcal V,
\]
(see \cite{klingenberg2012lectures}  \S1.1, Proposition~1.1.1). At $(x,v)\in T^1M$, this identifies
\[
T_{(x,v)}T^1M\cong T_xM\oplus v^\perp.
\]
Concretely, a tangent vector represented by a curve $s\mapsto(x(s),v(s))$ corresponds to
\[
(X,Y)=\left(\dot x(0),\left.D_sv(s)\right|_{s=0}\right).
\]
Here $X$ describes the motion of the base point, while $Y$ describes the covariant change of the unit vector. Since $|v(s)|=1$, we have $Y\perp v$.
Furthermore,
\[
\bar\alpha_{(x,v)}(X,Y)=\langle v,X\rangle.
\]
Restricting to $\xi$ therefore also imposes $X\perp v$. Along $\Gamma$, we obtain
\[
\xi_{\Gamma(t)}
\cong\dot\gamma(t)^\perp\oplus\dot\gamma(t)^\perp.
\]
In this splitting,
\begin{equation}\label{eqtwocomponentsymplecticsplitting}
\D\bar\alpha\bigl((X,Y),(\widehat X,\widehat Y)\bigr)
=\langle Y,\widehat X\rangle-\langle X,\widehat Y\rangle.
\end{equation}
Vertical vectors have the form $(0,Y)$, so their mutual symplectic pairings vanish.\\

\noindent The differential of the geodesic flow is
\begin{equation}\label{eqtwocomponentjacobiflow}
(\D\phi_t)_{\Gamma(0)}(X,Y)
=\bigl(J_{X,Y}(t),D_tJ_{X,Y}(t)\bigr),
\end{equation}
where $J_{X,Y}$ is the Jacobi field satisfying
\[
J_{X,Y}(0)=X,\qquad D_tJ_{X,Y}(0)=Y.
\]
This is the Jacobi-field description in Klingenberg (see Lemma~3.1.6. in \cite{klingenberg2012lectures}). The closed-geodesic index agrees here with the fixed-endpoint index. Indeed, $\D\phi_L=\mathrm{id}$ implies that every normal Jacobi field satisfies
\[
J(L)=J(0),\qquad D_tJ(L)=D_tJ(0).
\]
If $I$ denotes the periodic index form, integration by parts gives
\[
I(J,W)=\bigl[\langle D_tJ,W\rangle\bigr]_0^L=0
\]
for every periodic normal vector field $W$.
Given such a field $W$, choose the normal Jacobi field $J$ with $J(0)=W(0)$ and $D_tJ(0)=0$. Then $W-J$ vanishes at both endpoints. Consequently, the space of periodic normal variations is the direct sum of the fixed-endpoint variations and the Jacobi fields with $D_tJ(0)=0$; the latter summand lies in the kernel of $I$. The two index forms therefore have the same index. Since tangential variations contribute a nonnegative form, the fixed-endpoint Morse index theorem yields
\begin{equation}\label{eqtwocomponentconjugatecount}
\ind_{\Omega}=\ind(K_{i_0})=2(q_1+q_2)-2,
\end{equation}
where $\ind_{\Omega}(\gamma)$ denotes the fixed-endpoint Morse index of $\gamma$, equal to the number of conjugate points in $0<t<L$ counted with multiplicity. Now
\[
\Lambda(t):=G(t)\IR^2
\]
represents the image of the initial vertical space under the differential of the flow. Its intersection with the vertical space at time $t$ corresponds precisely to normal Jacobi fields vanishing at $0$ and $t$. Thus
\begin{equation}\label{eqtwocomponentverticaljacobi}
\dim_{\IR}\bigl(\Lambda(t)\cap\IR^2\bigr)
=\dim\bigl\{J:\ J\text{ is Jacobi},\ J(0)=J(t)=0,\ J\perp\dot\gamma\bigr\}.
\end{equation}
Define
\[
S(t)=G(t)G(t)^T.
\]
This matrix is symmetric and unitary, and
\[
\det S(t)=(\det G(t))^2.
\]
Hence
\begin{equation}\label{eqtwocomponentdoubledegree}
\deg(\det S)=2w.
\end{equation}
For a real vector $x$, the equality $S(t)x=x$ is equivalent to $x\in G(t)\IR^2$. The complex $1$-eigenspace of $S(t)$ is the complexification of its real part. Consequently,
\begin{equation}\label{eqtwocomponenteigenmultiplicity}
\dim_{\IC}\ker(S(t)-I)
=\dim_{\IR}\bigl(\Lambda(t)\cap\IR^2\bigr).
\end{equation}

The essential point is that, for each eigenvalue in $\IS^1$, the derivative of its angle is positive whenever the eigenvalue equals $1$. Write an eigenvalue locally as
\[
\lambda_j(t)=e^{i\theta_j(t)}.
\]
Fix $t_*$ with $\lambda_j(t_*)=1$, and let $x=x_j\in\IR^2$ be a corresponding nonzero eigenvector. Choose $y\in\IR^2$ such that
\[
z(t):=G(t)y,\qquad z(t_*)=x.
\]
In the frame $e(t)$, this represents a Jacobi field:
\[
e(t)z(t)=\bigl(J_j(t),D_tJ_j(t)\bigr),
\qquad
J_j(0)=J_j(t_*)=0.
\]
Since $G(t)$ is unitary,
\[
\begin{aligned}
S(t)\overline{z(t)}
&=G(t)G(t)^T\overline{G(t)}\,y\\
&=G(t)y=z(t).
\end{aligned}
\]
Differentiating at $t_*$ gives
\[
\dot S(t_*)x+S(t_*)\overline{\dot z(t_*)}
=\dot z(t_*).
\]
Because $S(t_*)x=x$ and $S(t_*)$ is unitary, $x^*S(t_*)=x^*$. Multiplying by $x^*$, we obtain
\begin{equation}\label{eqtwocomponentmatrixderivative}
\begin{aligned}
x^*\dot S(t_*)x
&=x^*\dot z(t_*)-\overline{x^*\dot z(t_*)}\\
&=2i\,\Im\bigl(x^*\dot z(t_*)\bigr).
\end{aligned}
\end{equation}
Put $V=D_tJ_j(t_*)$. Since $J_j(t_*)=0$, the Jacobi equation gives $D_t^2J_j(t_*)=0$. Thus the vector and its covariant derivative in the connection-induced splitting are respectively
\[
(0,V)\qquad\text{and}\qquad(V,0).
\]
Using \eqref{eqtwocomponentsymplecticsplitting} and the fact that $e(t)$ is unitary,
\begin{equation}\label{eqtwocomponentpositivepairing}
\begin{aligned}
\Im\bigl(x^*\dot z(t_*)\bigr)
&=\D\bar\alpha\bigl((0,V),(V,0)\bigr)\\
&=|V|^2.
\end{aligned}
\end{equation}
Differentiating the moving frame adds a vertical vector. Its pairing with the vertical vector $(0,V)$ is zero, so it does not affect this equality.\\

\noindent On the other hand, differentiating the eigenvalue equation at $\lambda_j(t_*)=1$ yields
\[
x^*\dot S(t_*)x=i\theta_j'(t_*)|x|^2.
\]
Together with \eqref{eqtwocomponentmatrixderivative} and \eqref{eqtwocomponentpositivepairing}, this gives
\begin{equation}\label{eqtwocomponentpositiveangles}
\theta_j'(t_*)|x_j|^2
=2|D_tJ_j(t_*)|^2>0.
\end{equation}
The inequality is strict because a Jacobi field whose value and derivative both vanish at one time is identically zero. If $1$ has multiplicity two, then $S(t_*)=I$. Choose real eigenvectors of the real symmetric matrix $-i\dot S(t_*)$. Its eigenvalues are the derivatives of the two angles, and the same calculation shows that both are positive.\\
It remains to count the degree. Choose continuous arguments on $[0,L]$,
\[
\lambda_j(t)=e^{i\theta_j(t)},\qquad \theta_j(0)=0.
\]
Since
\[
\det S(t)=e^{i(\theta_1(t)+\theta_2(t))},
\]
the definition of degree gives
\[
\deg(\det S)
=\frac1{2\pi}\sum_{j=1}^2\bigl(\theta_j(L)-\theta_j(0)\bigr).
\]
Equivalently,
\begin{equation}\label{eqtwocomponentsigneddegree}
\deg(\det S)
=\sum_{j=1}^2
\ \sum_{\substack{0\leq t<L\\\lambda_j(t)=1}}
\sgn\bigl(\theta_j'(t)\bigr),
\end{equation}
where the derivatives are taken along local eigenvalue branches.
Indeed, $\lambda_j(t)=1$ means that $\theta_j(t)$ is a multiple of $2\pi$. Passing such a value in the increasing direction contributes $+1$, while passing it in the decreasing direction contributes $-1$. Their signed sum gives the total number of turns. By \eqref{eqtwocomponentpositiveangles}, every sign in \eqref{eqtwocomponentsigneddegree} is positive. An eigenvalue $1$ of multiplicity $m$ therefore contributes $m$. Using \eqref{eqtwocomponentverticaljacobi} and \eqref{eqtwocomponenteigenmultiplicity}, we obtain
\begin{equation}\label{eqtwocomponentdegreecount}
\begin{aligned}
\deg(\det S)
&=\sum_{\substack{0\leq t<L\\1\in\operatorname{spec}S(t)}}
\dim_{\IC}\ker(S(t)-I)\\
&=\ind_{\Omega}(\gamma)+2.
\end{aligned}
\end{equation}
The interior times contribute $\ind_{\Omega}$, while $S(0)=I$ contributes $2$. We exclude $L$, since $0$ and $L$ represent the same point of the time circle.
Thus \eqref{eqtwocomponentconjugatecount}, \eqref{eqtwocomponentdoubledegree}, and \eqref{eqtwocomponentdegreecount} give
\[
2w=\deg(\det S)=\ind_{\Omega}(\gamma)+2=2(q_1+q_2),
\]
and hence
\begin{equation}\label{eqtwocomponentwindingvalue}
w=q_1+q_2.
\end{equation}

\medskip
\noindent\textbf{Claim 3.}
The orbifold first Chern class equals $w$ times the Euler class. In the notation of \eqref{eqtwocomponentchernmultiplelcm}, this means
\[
\kappa=w.
\]

We construct a sphere in the classifying space on which the symplectic class evaluates to $1$ and the first Chern class evaluates to $w$.
Consider
\[
Q:T^1M\times_{\IS^1} E\IS^1\longrightarrow T^1M/\IS^1,
\qquad
Q([u,a])=\Pi_{\IS^1}(u).
\]
The disk
\[
\widetilde D_+:D_+\longrightarrow T^1M\times_{\IS^1} E\IS^1,
\qquad
\widetilde D_+(z)=[D(z),a_0],
\]
has boundary contained in the fiber over
\[
p=\Pi_{\IS^1}(\Gamma).
\]
Because $\Gamma$ is a principal orbit,
\[
Q^{-1}(p)\cong E\IS^1
\]
is contractible. We can therefore fill the boundary by a second disk $D_-$ contained in this fiber. Gluing the disks with opposite boundary orientations gives
\[
f_D:\IS^2\longrightarrow T^1M\times_{\IS^1} E\IS^1.
\]
Set
\[
\beta_D=(f_D)_*[\IS^2],
\qquad
E_D=f_D^*\mathscr E.
\]
The bundle $E_D$ is an ordinary complex vector bundle over $\IS^2$, and naturality gives
\begin{equation}\label{eqtwocomponentchernnaturality}
\left\langle c_1^{\mathrm{orb}}\bigl(T(T^1M/\IS^1)\bigr),\beta_D\right\rangle
=\langle c_1(E_D),[\IS^2]\rangle.
\end{equation}
Over $D_+$, the bundle $E_D$ has the frame $e$. Over $D_-$, it is identified with the constant vector space $T_p(T^1M/\IS^1)$. Along the common boundary, the corresponding frame upstairs is
\[
e_{\mathrm{flow}}(t)=(\D\phi_t)e(0),
\]
because
\[
\D\Pi_{\IS^1}\circ\D\phi_t=\D\Pi_{\IS^1}.
\]
The definition of $G$ gives
\[
e_{\mathrm{flow}}(t)=e(t)G(t).
\]
Taking complex wedges,
\begin{equation}\label{eqtwocomponentdettransition}
s_{\mathrm{flow}}(t)=\det G(t)\,s_D(t),
\qquad
s_D=e_1\wedge_{\IC}e_2.
\end{equation}
Extend the nonzero determinant section on $D_-$ across $D_+$, allowing isolated zeros in the interior. In the disk trivialization, write the extension as
\[
\sigma=f\,s_D.
\]
Its boundary value is
\[
f|_{\partial D_+}=\det G.
\]
The signed number of zeros of a section of a complex line bundle is its first Chern number. To make the signs explicit, let $B_\ell$ be small positively oriented disks around the zeros. The local contribution of a zero is
\[
\deg\left(\left.\frac{f}{|f|}\right|_{\partial B_\ell}\right).
\]
A zero modeled by $f(z)=z$ contributes $+1$, while one modeled by $f(z)=\bar z$ contributes $-1$.
Since $c_1(E_D)=c_1(\det_{\IC}E_D)$, we obtain
\begin{equation}\label{eqtwocomponentchernzerocount}
\begin{aligned}
\langle c_1(E_D),[\IS^2]\rangle
&=\sum_\ell\deg\left(\left.\frac{f}{|f|}\right|_{\partial B_\ell}\right)\\
&=\deg\left(\left.\frac{f}{|f|}\right|_{\partial D_+}\right)\\
&=\deg(\det G)=w.
\end{aligned}
\end{equation}
The second equality follows by removing the small disks: the degree on the outer boundary equals the sum of the degrees on their positively oriented boundaries. The disc $D_-$ projects to a point in the quotient and therefore contributes zero symplectic area. Stokes' theorem gives
\begin{equation}\label{eqtwocomponentstokeslcm}
\begin{aligned}
\langle[\omega_{orb}],\beta_D\rangle
&=\frac1L\int_{D_+}D^*\D\bar\alpha\\
&=\frac1L\int_\Gamma\bar\alpha\\
&=1,
\end{aligned}
\end{equation}
because $\bar\alpha(\dot\Gamma)=1$ and $\Gamma$ has period $L$.
Evaluating \eqref{eqtwocomponentchernmultiplelcm} on $\beta_D$ and using \eqref{eqtwocomponentchernnaturality}, \eqref{eqtwocomponentchernzerocount}, and \eqref{eqtwocomponentstokeslcm} proves $\kappa=w$. Together with Claim~2, this gives $\kappa=q_1+q_2$ and therefore
\begin{equation}\label{eqtwocomponentchernidentitylcm}
c_1^{\mathrm{orb}}\bigl(T(T^1M/\IS^1)\bigr)
=(q_1+q_2)[\omega_{orb}].
\end{equation}

\medskip
\noindent\textbf{Claim 4.}
The exceptional components give
\[
\kappa=\frac{2(q_1+q_2)}
{\#\pi_1(C_1)+\#\pi_1(C_2)}.
\]
By Lemma~\ref{lembasisexistence}, the two singular surfaces have classes
\[
[C_i/\IS^1]=E_1+l_iE_2,
\qquad
E_1^2=E_2^2=0,
\qquad
E_1\cdot E_2=1.
\]
Since the surfaces are disjoint,
\[
0=[C_1/\IS^1]\cdot[C_2/\IS^1]=l_1+l_2.
\]
Their self-intersections therefore satisfy
\begin{equation}\label{eqtwocomponentselfintersectionsum}
[C_1/\IS^1]^2+[C_2/\IS^1]^2=2(l_1+l_2)=0.
\end{equation}
We continue to write $C_i/\IS^1$ for the corresponding embedded symplectic sphere in the smooth symplectic resolution.

The codimension-two canonical-class formula of Chen \cite[Proposition~3.2]{chen2018resolving}, together with the adjunction formula \cite[Theorem~2.6.4, p.~36]{McDuffSalamon2012}, yields
\begin{equation}\label{eqtwocomponentorbifoldadjunction}
\left\langle
c_1^{\mathrm{orb}}\bigl(T(T^1M/\IS^1)\bigr),[C_i/\IS^1]
\right\rangle
=2+\frac{[C_i/\IS^1]^2}{q_i}.
\end{equation}
Indeed, adjunction gives $2+[C_i/\IS^1]^2$ for the embedded sphere in the smooth resolution, and Chen's formula subtracts
\[
\left(1-\frac1{q_i}\right)[C_i/\IS^1]^2.
\]
This gives \eqref{eqtwocomponentorbifoldadjunction}. On the other hand, the circle bundle
\[
\IS^1/(\IZ/q_i\IZ)\longrightarrow C_i\longrightarrow C_i/\IS^1
\]
has Euler number $\#\pi_1(C_i)$ with respect to the orientation induced by $\omega_{orb}$. Its normalized connection is $\frac{q_i}{L}\bar\alpha$, which gives
\begin{equation}\label{eqtwocomponenteulerlcm}
q_i\left\langle[\omega_{orb}],[C_i/\IS^1]\right\rangle
=\#\pi_1(C_i).
\end{equation}
Combining \eqref{eqtwocomponentchernmultiplelcm}, \eqref{eqtwocomponentorbifoldadjunction}, and \eqref{eqtwocomponenteulerlcm}, we obtain
\[
\kappa\,\#\pi_1(C_i)=2q_i+[C_i/\IS^1]^2,
\qquad i=1,2.
\]
Adding these identities and using \eqref{eqtwocomponentselfintersectionsum} yields
\[
\kappa\bigl(\#\pi_1(C_1)+\#\pi_1(C_2)\bigr)=2(q_1+q_2),
\]
which proves Claim~4.

\medskip
\noindent Finally, comparing Claims~3 and~4, we have
\[
q_1+q_2
=\frac{2(q_1+q_2)}
{\#\pi_1(C_1)+\#\pi_1(C_2)}.
\]
Since $q_1+q_2>0$, it follows that
\[
\#\pi_1(C_1)+\#\pi_1(C_2)=2.
\]
This is impossible, since both $\#\pi_1(C_1)$ and $\#\pi_1(C_2)$ are positive even integers. Thus the case $k=2$ is also excluded. The non-trivial isotropy set must therefore be empty, so the geodesic-flow action is free. Equivalently, all prime geodesics have the same length. This completes the proof of Theorem~\ref{bergconj}.
\end{proof}